\documentclass[12pt]{article}
\title{Test vectors for some ramified representations}
\author{Vinayak Vatsal}
\date{\today}
\usepackage{commands}
\begin{document}
\maketitle

\begin{abstract}We give an explicit construction of test vectors for $T$-equivariant linear functionals on representations $\Pi$  of $GL_2$ of a $p$-adic
field $F$, where $T$ is a non-split
torus. Of particular interest is the case when both the representations are ramified; 
we completely solve this problem for principal series and Steinberg representations
of $GL_2$, as well as for depth zero supercuspidals over $F=\qq_p$. A key ingredient is a theorem of Casselman and Silberger, 
which allows us to quickly reduce almost all cases
to that of the principal series, which can be analyzed directly. 
Our method shows that the only genuinely difficult cases are the characters of $T$ which occur in the primitive part 
(or ``type") of $\Pi$ when $\Pi$ is supercuspidal.
The method to handle the depth zero case is 
based on modular representation theory, motivated by considerations from Deligne-Lusztig theory and the de Rham cohomology of
Deligne-Lusztig-Drinfeld curves. The proof also reveals some interesting features related to the Langlands correspondence in characteristic 
$p$. We show in particular that the test vector problem has an obstruction in characteristic $p$ beyond the root number criterion of 
Waldspurger and Tunnell, and that it exhibits an unexpected dichotomy related
to the weights in Serre's conjecture and the signs of standard Gauss sums. 
\end{abstract}

\pagebreak
\tableofcontents

\section{Introduction}

Let $F$ be a finite extension of $\qq_p$, for a rational prime $p$.
Let $G$ denote the algebraic group $GL_2$ over $F$  and let
$T\subset GL_2$ denote a maximal torus. 
Let $\chi$ denote a $1$-dimensional character of $T(F)$. Let $\Pi$ denote an
irreducible smooth admissible infinite dimensional representation of $G(F)$, with central 
character $\omega$. Suppose that $\chi$ agrees
with $\omega$ on $F^*\subset T(F)$. The space $\hom_T(\Pi, \chi)$ has dimension
$1$ when it non-empty, which we assume to be the case. Let $\ell_\chi$ be a nonzero
element of $\hom_T(\Pi, \chi)$. Then a \emph{test vector}
for $\ell_\chi$ is a vector $v\in \Pi$ such that $\ell_\chi(v)$ is nonzero. The study of 
test vectors was initiated by Gross and Prasad in \cite{grp91}. The main
results in \cite{grp91} treat the cases where the ramification
of $\Pi$ and $\chi$ is disjoint. It was already well-known, even before Gross-Prasad, 
that the case where $\Pi$  and $\chi$ are both ramified is significantly more complicated. 
Indeed, this problem comes up in producing explicit formulae for the local Rankin-Selberg integrals
of a tensor product $\Pi_1\otimes\Pi_2$ when $p$ is a bad prime for both $\Pi_1$ and $\Pi_2$. Many of
these ramified cases were treated in \cite{fmp}, \cite{cst}, and \cite{yueke-hu}, but the general problem remains open. 

From the point of view of number theory, having a recipe to produce a test vector at bad primes
is highly desirable, in order to give completely explict formulae of Gross-Zagier type.
The consideration of bad places is unavoidable in Iwasawa theory and $p$-adic deformation theory, but 
many results of Gross-Zagier type impose unpleasant restriction on the level, or 
are proven only up to a nonzero constant, arising from unknown local integrals at the bad places. 
Furthermore, applications to number theory demand that the test vector 
arising should be related in some elementary way to the newform theory of Atkin-Lehner. As expressed
in \cite{fmp}, one would like to exhibit a test vector in the form of an explicit translate of the Atkin-Lehner
new vector, and one would like further to have a translate whose definition depends only on elementary
 ramfiication data such as the conductor and central character of $\Pi$ and the conductor of the character $\chi$.
In the original setting of \cite{grp91}, the emphasis was placed on finding a test vector that was fixed 
under the group of units of a suitably chosen Eichler order. 

The original motivation for this paper was to produce test vectors in some cases that were needed in the construction
of \cite{vat17a}, since these results seemed unavailable
in the literature. In the process of filling in the missing cases, we discovered that it was possible to give substantially simpler proofs
of some of the results of \cite{fmp}, and also give a number of new results that were not covered in \cite{fmp} at all. 
The results of \cite{fmp} are comprehensive in the case of a split torus $T$; our results in this paper pertain to the case where the torus
$T$ corresponds to a quadratic field extension $K/F$. For the remainder of this paper, we assume that $T=K^*$, for a quadratic
field extension $K/F$. 

Our first result is quite general. We exhibit, in full generality, a test vector $v_\chi$ such that $\ell_\chi(v_\chi)\neq 0$, provided that
the ramification of $\chi$ is large compared to that of $\Pi$. The test vector we supply is an obvious translate of the new vector, and the definition
depends only on elementary ramification data. This result is already contained in \cite{fmp}, but our proof is simpler, and
gives a pleasant explicit bound on the ramifications in question. The key new ingredient is systematic use of a theorem of Silberger \cite{silberger}
and Casselman \cite{cass-restriction}.
Casselman and Silberger show that the restriction of $\Pi$ to $GL_2(\oo_F)$ decomposes as the sum of a finite dimensional piece (the primitive piece, or ``type"),
which we denote as $V(\Pi)$,
plus an infinite dimensional piece, denoted $V(\omega)$, that depends only on the central character $\omega$. The torus $T(F)$ is compact mod centre, under our 
assumptions, and one obtains a corresponding decomposition of $\Pi$ into characters of $T$. Each character appears with multiplicity one 
so we can classify the characters of $T$ which appear according to whether they appear in $V(\Pi)$ or $V(\omega)$. 

Then, our first result solves the test vector problem for those characters of $T$ which appear in $V(\omega)$. The constituents of $V(\omega)$
as a representation of $GL_2(\oo_F)$ are essentially induced, and producing test vectors is an elementary problem. Since $V(\omega)$ depends
only on $\omega$, our results here are valid for all $\Pi$. One way of interpreting this result is that the decomposition
 $V=V(\Pi)\oplus V(\omega)$ is a manifestation of the phenomenon that the epsilon factors of twists of $\Pi$ stabilize once the twist is sufficiently ramified; 
 since the existence of a nonzero $\chi$-equivariant functional $\ell_\chi$ is governed by the value of a certain epsilon factor, one sees that the criterion
 which governs the existence of $\ell_\chi$ is automatic once $\chi$ is sufficiently ramified, and the existence of a good test vector is more or less automatic
 as well. Since the results depend only on $\omega$, one has only to analyze the principal series, and the results are given in Section \ref{pseries} below. 
 Another perspective was pointed out by Prasad -- there is only one singular point for the group $PGL_2$ which is the identity,
and if ones subtracts the characters of two representations, one gets a function which has no singularity and is locally 
constant everywhere, and whose restriction to any compact open subgroup such as $GL_2(\oo)$ 
is the character of a finite dimensional representation.
 
 The more interesting case is that of characters $\chi$ which occur in the type $V(\Pi)$. Here the root number condition controlling the existence of 
 $\ell_\chi$ is non-empty, and not easy to make explicit. Any construction of the test vector has to take in to account the root numbers, and cannot
 simply depend on the levels of the characters and representations. This kind of analysis 
 is not too hard for the case of special or principal series representations, where the type $V(\Pi)$ is
 rather simple, but it is highly nontrivial in the case of supercuspidals, where the finite dimensional representation $V(\Pi)$ is inflated from an irreducible
 representation of the finite group $GL_2(\oo_F/\varpi^n\oo_F)$, for suitable $n$. Here $\varpi$ is a uniformizer of $\oo_F$. 
 
 The most interesting results of this paper give a detailed analysis of the test vector problem for 
 characters $\chi$ of an unramified quadratic extension
 $K/\qq_p$ and 
 the case of depth zero supercuspidals  $\Pi$. This means that
$V(\Pi)$ has dimension $p-1$, and is inflated from an irreducible cuspidal representation of the $GL_2(\ff_p)$. 
  To give an idea of how basic the test vector problem is, consider the case that the central character $\omega$ is trivial; then 
 the trivial character of $T$ occurs in $V(\Pi)$, as shown by Gross \cite{gro88}, but the results of this 
 paper are the first to give an natural test vector for the unique
 $T$-equivariant linear functional on $\Pi$.  We remark here that this case is briefly treated in \cite{grp91} and \cite{gro88} but the test
 vector there is not related in any way to the invariants of an Eichler order, and is not connected in any way to the newforms of number theory.
 Our results are quite surprising, 
 and are intimately related to representation theory in characteristic $p$.
 
To state the results we recall some background material on representations of $GL_2$ over a finite field, as 
given in the book \cite{ps-book}. It is shown in that book
that the cuspidal complex representations of $GL_2(\ff_q)$ are parametrized by characters $\nu:\ff_{q^2}^*\rightarrow\cc^*$
such that $\nu\neq \nu^q$. Any such $\nu$ gives rise to a $q-1$-dimensional representation 
$\rho(\nu):GL_2(\ff_q)\rightarrow GL_{q-1}(\cc)$ with central character
$\ff_q^*\rightarrow\cc^*$ given by the restriction of $\nu$ to $\ff_q^*$. 
We have $\rho(\nu)\cong\rho(\nu^q)$, and this is the only possible isomorphism
between the various representations $\rho(\nu)$, for varying $\nu$. 

Now let $K/F$ denote the unramified quadratic extension. Since $p$ is odd, we may and do identify $K$ with the subgroup of matrices $\begin{pmatrix}
a & b\\ b\xi & a\end{pmatrix}$ with $a, b\in F, a^2-\xi b^2\neq 0$ where $\xi\in F^*$ is a non-square unit. 
Reduction modulo $\varpi$ gives rise to an embedding
$\ff_{q^2}^*\subset GL_2(\ff_q)$. Let $\nu$ denote a character of $\ff_{q^2}^*$ as above, and let $\omega$ be any character
of $F^*$ which agrees with (the inflation of) $\nu$ on $\oo_F^*$. We remark that $\nu$ is determined
 on $\oo_K^\times$ by inflation, but that we can extend $\nu$ to $K^*$ by requiring it to agree with $\omega$
 on $F^*$.  
Then there exists an irreducible smooth admissible
representation $\Pi$ of $GL_2(F)$ with central character $\omega$, such that the type $V(\Pi)$ is isomorphic to
$\rho(\nu)$. 

Consider a character $\chi$  of $T=K^*$ such that $\chi$ agrees with 
$\omega$ on $F^*$.  Assume further that $\chi$ is trivial on the subgroup $1+\varpi_K\oo_K$ of principal units in $K$. 
Then we may identify $\chi$ with a character of $\ff_{q^2}^*$ which agrees with $\nu$ on $\ff_q^*$ (the central character
being fixed). 
It is a basic fact about cuspidal representations of 
$GL_2(\ff_q)$ that each such character $\chi\neq \nu,\nu^q$ occurs with multiplicity one in the representation $\rho(\nu)$. 
Thus, for each choice
of $\chi$, we get a distinguished line in the space $V(\Pi)$. As was observed by Gross and Prasad, 
a vector $v$ in $V$ is a test vector for $\chi$
if and only if it has nonzero projection on this distinguished line (here projection means orthogonal projection 
with respect to the action of the compact group $\oo_K^*$). 

Now let $v^\text{new}\in\Pi$ denote the new vector provided by Casselman's local version of Atkin-Lehner theory \cite{cass73}. By definition,
$v^\text{new}$ is a nonzero vector on the unique line in $V$ such that the subgroup $\begin{pmatrix}  a& b \\ \varpi^2 c & d\end{pmatrix}$,
with $a, d\in \oo_F^*$, $b, c\in\oo_F$, acts via the character $\omega(a)$. 
A basic observation is that $v^\text{new}$ lies in $V(\omega)$
and has projection zero to $V(\Pi)$. Thus $v^\text{new}$ is \emph{never} a test vector for \emph{any} character 
$\chi$ appearing in $V(\Pi)$, and to produce
a plausible candidate for test vector, one has to translate $v^\text{new}$ to put it inside $V(\Pi)$. Thus, consider the vector 
$$v_1=\begin{pmatrix} \varpi^{-1}& 0 \\ 0 & 1\end{pmatrix}\cdot v^\text{new}$$
where the action of the matrix is given by the representation $\Pi$. An easy calculation show that $v_1$ is invariant under the group 
of matrices congruent to the identity modulo $\varpi$. Thus $v_1\in V(\Pi)$ and in fact $v_1$ is
an eigenvector for the subgroup $D^0$ of matrices of the form $\begin{pmatrix} a & 0 \\ 0 & 1\end{pmatrix}$, $a \in \oo_F^*$, 
with eigenvalues given by $\omega(a)$. Thus, to determine whether  $v_1$ 
is a test vector for $\chi$, one has to determine
whether an eigenvector for the \emph{split} diagonal torus has nonzero projection on an eigenspace of a \emph{non-split} torus. 
It is this tension between the split and non-split
tori that makes the problem hard. 

More generally, one has to solve the following problem in finite group theory. Let $\nu:\ff_{q^2}^*\rightarrow\cc^*$ be a character
such that $\nu\neq \nu^q$, and let $\rho(\nu)$ denote corresponding cuspidal representation of $GL_2(\ff_q)$. Let $\mu$ denote
a character of $\ff_q^*$, and let $v_\mu$ denote an eigenvector for the group $D^0=\begin{pmatrix} a & 0 \\ 0 & 1\end{pmatrix}$,
$a\in \ff_q^*$, with eigenvalues under $\rho(\nu)$ given by
$$\begin{pmatrix} a & 0 \\ 0 & 1\end{pmatrix}\cdot v_\mu= \mu(a) v_\mu.$$
Such an eigenvector always exists, for any $\mu$. In fact, $\rho(\nu)$ decomposes as the sum of $q-1$ invariant lines, 
each being an eigenspace for a unique character $\mu$.
Let $\chi:\ff_{q^2}^*\rightarrow \cc^*$ denote a character of the nonsplit 
torus $T$ such that $\chi$ agrees with $\nu$ on $\ff_q^*$, such that $\chi\neq\nu, \nu^q$. Here we identify the nonsplit
torus with $2\times 2$ matrices as explained above. Then, is the projection of $v_\mu$ to the
$\chi$-eigenspace of $T$ nonzero? It is this problem that we solve in this paper, under the key hypothesis that $q=p$, 
so that we are dealing with the prime field. While our methods
give some information for general $q$, the combinatorics are quite complicated, and we have not pursued a general theorem. 

To give the reader a taste of our results, we give a sample in an illuminating special case, which arises from elliptic curves with
conductor $p^2$. Suppose that $\Pi$ is a depth zero supercuspidal representation of $GL_2(\qq_p)$ with trivial central character 
$\omega$. Then one
can ask whether or not the vector $v_\text{triv}$ corresponding to the trivial character $\mu=1$ is a test vector for some given $\chi$,
for instance, for $\chi=1$. Another obvious candidate for a test vector is the vector $v_\text{quad}$ corresponding to the choice of
of $\mu=\mu^\text{quad}$, the nontrivial quadratic character of $\ff_p^*$. The somewhat surprising 
answer is that both choices $v_\text{triv}$ and $v_\text{quad}$ work, but only half the time. 

To state our result precisely, we need some way to label the characters of $T$ which arise. Thus 
fix any $\chi$ such that $\chi$ is trivial on $\qq_p^*$ and on the principal units $1+p\oo_K$. Any such
character is defined on $\oo_K^*/1+p\oo_K\cong\ff_{p^2}^*$ and takes values in the group of $\mu_{p^2-1}$ of complex $p^2-1$-th 
roots of unity. If we fix a prime 
$\pp$ of $\qbar$ above $p$, we may identify 
$\chi$ with a character $\ff_{p^2}^*\rightarrow\ff_{p^2}^*$. Then $\chi$ can be expressed in the form $x\mapsto x^{a+pb}$,
where $0\leq, a, b\leq p-1$. This expression is unique except for the case $(a,b)=(0,0)$ and $(a, b)=(p-1, p-1)$. We fix 
the former choice in that case. We write $\chi=\chi(a, b)$. We remind the reader that the integers $(a, b)$ 
depend on a choice of $\pp$. Given a representation $\rho(\nu)$ attached to a 
character $\nu$ such that $\nu\neq\nu^p$, the integers $a, b$ determined by $\nu$ must satisfy
 $a\neq b$, and so we may assume, by switching $\nu$ and $\nu^p$ if
neccessary, that $a > b$. If $\nu$ is trivial on $\ff_p^*$, then we must have $a+b\equiv 0\pmod{p-1}$, namely, $a=p-1-b$. 
Note that $a=b=0$ is excluded since $\nu\neq\nu^p$. 

Now let $\chi$ denote any character of $K^*/\qq_p^*(1+p\oo_K)\cong\ff_{p^2}^*/\ff_p^*$, 
and write $\chi=\chi(r, s)$ as above, with $0\leq r, s, \leq p-1$ as above. 
Of course, for general
$\chi$ we may have $r=s$ or $r < s$. However, we still have $r+s\equiv0\pmod{p-1}$, which means $r+s=p-1$ unless
$r=s=0$. In particular, $\chi$ is determined by the value of $r$ except for $(0, 0)$ and $(0, p-1)$.
Thus there are $p+1$ choices for $\chi$, as expected, including $\chi=\nu, \nu^p$. Two of the $p+1$
possible characters $\chi$ factor through the norm, namely, the trivial character and the unique nontrivial character
of order $2$. The remaining $p-1$ characters are primitive in the sense that they do not factor through the norm.

To get a better label for $\chi(r, s)$ define $ k_\chi = r + sp$, so that $\chi(x)=x^{k_\chi}$. The integer $k_\chi$ is closely related to
 the weight arising in Serre's conjecture, but we will not need this connection here. Note that $k_\nu = a + bp$,
 and $k_{\nu^p}=b + ap$. Since $\nu\neq\nu^p$, we must have $k_\nu\neq k_{\nu^p}$; in fact, one has $k_{\nu^p}-k _\nu
 =(a-b)(p-1) > 0$, so $k_{\nu^p} > k_\nu$. 
 Furthermore, if $\rho(\nu)$ has trivial central 
 character, then $\nu$ is trivial on $\ff_p^*$, so that $k_\nu= bp+a\equiv 0 \pmod{p-1}$. Thus the integers $k_\chi$ run
 through the integers $0, p-1, 2(p-1), \dots, p(p-1)$. Since $a+b=p-1$, one has 
 $k_\nu=a+bp = p-1-b +bp=(p-1)(b+1)$ and $k_{\nu^p}=b+ap = (p-1)(a+1)$.  
 We label $\chi$ with  the integer $k_\chi$, and we shall say that $\chi$ has \emph{Type 1} if $$k_\nu < k_\chi < k_{\nu^p}.$$
Otherwise we say $\chi$ is of Type 2. Observe that the trivial character $\chi$ has $k_\chi=0$ and has Type 2, while the nontrivial 
quadratic character has $k_\chi=(p^2-1)/2$ and has Type 1.

Write $\psi$ to denote a nontrivial additive character of $\ff_p$. Then if $\chi$ is as above, let $G(\chi)$ denote
the standard finite field Gauss sum with respect to $\psi$ for $\chi$ (viewed as a complex character of $\ff_{p^2}^*$) so that $G(\chi)=
\sum_{u\in\ff_{p^2}}\chi(u)\psi(\tr(u))$. Then one has the following lemma, which will be proved in Section \ref{supercusp}.

\begin{lmm} Suppose $\chi$ does not factor through the norm, and that $\chi$ is trivial on $\ff_p^*$. 
Then $\epsilon_p(\chi) :=\frac{G(\chi)}{p}= \chi(\sqrt{\xi})=\sqrt{\xi}^{k_\chi}$, where $\sqrt{\xi}\in\ff_{p^2}^*$ is any element of trace
zero whose square is a non-square element $\xi$ in $\ff_p^*$. 
\end{lmm}

With this lemma in hand, we extend the definition of $\epsilon_p(\chi)$ to all $\chi$ by setting $\epsilon_p(\chi)=\chi(\sqrt{\xi})$ when
$\chi$ is either trivial or quadratic. Then our result is as follows.

\begin{thm}
\label{gauss-sum-intro} Let $\Pi$ denote a depth zero supercuspidal representation of $GL_2(\qq_p)$ associated to a character
$\nu:\ff_{p^2}^*\rightarrow\cc^*$. 
Assume that $\Pi$ has trivial central character, so the type $\rho(\nu)$ is such that $\nu=\nu(a, b)$ with $a > b$ and $a+b=p-1$. 
Let $\chi$
denote a character of  $K^*/\qq_p^*(1+p\oo_K)$, identified with a character of $\ff_{p^2}^*/\ff_p^*$. Then
 the following statements hold.
\begin{enumerate}
\item The vector $v_\text{triv}$ is a test vector for the trivial character $\chi=1$ if and only if $\epsilon_p(\chi)\neq\epsilon_p(\nu)$.
The character $\chi=1$ is of Type 2.
\item The vector $v_\text{quad}$ is a test vector for the nontrivial quadratic character $\chi$ of $K^*/\qq_p^*(1+p\oo_K)$  
if and only if $\epsilon_p(\chi)\neq\epsilon_p(\nu)$. The nontrivial quadratic character $\chi$ is of Type 1. 
\item More generally, the vector $v_\text{triv}$ is a test vector for all Type 2 characters $\chi$ such that $\epsilon_p(\chi)\neq \epsilon_p(\nu)$
and the vector $v_\text{quad}$ is a test vector for all Type 1 characters $\chi$ such that $\epsilon_p(\chi)\neq \epsilon_p(\nu)$.
\item Neither $v_\text{triv}$ nor $v_\text{quad}$ is a test vector for any character $\chi$ with  $\epsilon_p(\chi) = \epsilon_p(\nu)$.
\end{enumerate}
\end{thm}

One can give a similar result for general depth zero supercuspidals when $F=\qq_p$, as follows. Let $\Pi$ denote a depth
zero supercuspidal representation of $GL_2(\qq_p)$ with central character $\omega$, associated to the type $\rho(\nu)$.
Write $\nu=\nu(a, b)$ as above. Here we understand that $\nu$ is a character of $\ff_{p^2}^*$, extended to $K^*$ by 
virtue of the central character $\omega$, and the integers $a, b$ depend on a choice of prime above $p$ in $\qbar$,
and an identification of the complex $p^2-1$ roots of unity with $\ff_{p^2}^*$.
In particular, if $x\in\ff_p^*$, then we have $\nu(x)=\omega(x)=x^{a+pb}=x^{a+b}$, by definition. Given a character $\chi$
of $\oo_K^*/(1+p\oo_K)=\ff_{p^2}^*$ we extend it to $K^*$, since compatibility with $\omega$ specifies the value $\chi(p)=\omega(p)$. 
It will be notationally convenient to deal only with characters of the finite groups $\ff_p^*$ and $\ff_{p^2}^*$, and we generally 
will do this without comment in the sequel, since $\omega$ will be fixed.

Once again, we have $p+1$ characters $\chi$, including $\nu, \nu^p$ of $\ff_{p^2}$
which agree with $\nu$ on $\ff_p^*$.
We extend $\chi$ to characters of $K^*$ as above. Write $\chi=\chi(r, s)$, with $0\leq, r, s \leq p-1$ as explained
above, and set $k_\chi=r+sp$, then $k_\chi\equiv k_\nu\pmod{p-1}$, and we may write $k_\chi=k_\nu + t_\chi(p-1)$
for some $t_\chi$ which is uniquely determined modulo $p+1$. Once again we can classify characters of Type 1 and Type
2, with the exactly the same definition as before, namely, by requiring that the characters of Type 1 are those whose
weights lie between the weights of $\nu$ and $\nu^p$. 

Assume first that $\omega$ is even, in the sense that $\omega(-1)=1$, or -- equivalently -- that $a+b$ is even. Let
$\mu_1$ and $\mu_2$ denote the two possible characters of $\ff_p^*$ such that $\mu_i^2=\omega$, and let
$\chi_i=\textbf{N}\cdot\mu_i$ be the composition with the norm. For concreteness, we take $\mu_1(x) = x^{(a+b)/2}$
and $\mu_2(x) = x^{(a+b+p-1)/2}$. 
Let $\mu_3$ denote the character of $\ff_p^*$ defined
by $x\mapsto x^{(a+b)/2-1}$ and let $\mu_4(x)=x^{(a+b+p-1)/2-1}$.

\begin{thm}
\label{even-case-intro} Let $\Pi$ denote the depth zero supercuspidal representation associated to a character
$\nu:\ff_{p^2}^*\rightarrow\cc^*$. Let $\nu=\nu(a, b)$ with $a+b$ even, $a>b$. Let $\chi$ be a character of the unramfied
quadratic extension $K/\qq_p$ which agrees with $\nu$ on $\ff_{p^2}^*$. 
 Then the following statements hold:
\begin{enumerate}
\item The vector $v_{\mu_1}$ is a test vector for the character $\chi_1$ if and only if $t_{\chi_1}$ is odd. The character
$\mu_1$ is of Type 1.
\item The vector $v_{\mu_2}$ is a test vector for the character $\chi_2$ if and only if $t_{\chi_2}$ is odd. The character
$\mu_2$ is of Type 2. 
\item  More generally, the vector $v_{\mu_1}$ is a test vector for all characters $\chi$  of Type 1 such that $t_\chi$ is odd, and
the vector $v_{\mu_2}$ is a test vector for all characters $\chi$ of Type 2  such that $t_\chi$ is odd. 
\item Neither $v_{\mu_1}$ nor $v_{\mu_2}$ is a test vector for any character $\chi$ with  $t_\chi$ even.
\item The vector $v_{\mu_3}$ is a test vector for all characters $\chi$ of Type 1 such that $t_\chi$ is even.
\item The vector $v_{\mu_4}$ is a test vector for all characters $\chi$ of Type 2 such that $t_\chi$ is even.

\end{enumerate}
\end{thm}

We remark that the statements in the theorem above degenerate when $p=3$, and $\rho(\nu)$ has dimension 2, in the 
sense that some of the statements are empty. There are only two candidates for $\chi$ in this case, namely, the trivial character and
the nontrivial quadratic.

Finally, we give the theorem in the odd case. Let $\nu:\ff_{p^2}^*\rightarrow\cc^*$ be a character such that $\nu=\nu(a, b)$ with $a+b$ 
odd and $\nu\neq\nu^p$. Let $\omega$ be a character of $\qq_p^*$ such that the composition of $\omega$ with the 
Teichm\"uller lift of  $\ff_p^*$ agrees with $\nu$. Let $\Pi$ denote the depth zero supercuspidal associated to $\nu$ 
and the choice of $\omega$. Let $\chi$ denote a character of the unramfied quadratic extension $K$ of $\qq_p$ 
which agrees with $\omega$ on $\qq_p^*$. We can identify $\chi$ with a character of $\ff_{p^2}^*$ as above,
and define the notion of Type 1 and Type 2, exactly as before. 

Let $\mu_1, \mu_2$ be characters of $\ff_p^*$ defined by $\mu_1(x)=x^{(a+b-1)/2}$ and $\mu_2(x)x^{(a+b+p)/2}$,
and let $v_{\mu_i}$ denote corresponding eigenvectors in the type $\rho(\nu)$ of $\Pi$. 

\begin{thm}
\label{odd-case-intro} Let $\Pi$ denote the depth zero supercuspidal representation associated to a character
$\nu:\ff_{p^2}^*\rightarrow\cc^*$. Let $(a, b)$ denote the pair associated to $\nu$ and assume that $a+b$ is odd, $a > b$. 
Let $\chi$ be a character of $K^*$ which agrees with $\nu$ on $\qq_p^*$ and which is distinct from $\nu, \nu^p$
on $\ff_{p^2}^*$. 
 Then the following statements hold if $a-b > 1$:
\begin{enumerate}
\item The vector $v_{\mu_1}$ is a test vector for $\chi$ if $\chi$ has Type 1.
\item The vector $v_{\mu_2}$ is a test vector for $\chi$ if $\chi$ has Type 2.
\end{enumerate}
If $a-b=1$ then the vector $v_{\mu_2}$ is a test vector for all $\chi$.
\end{thm}

\subsection*{Deligne-Lusztig theory, Serre weights, and the reduction to characteristic $p$}
We conclude the introduction with some comments on the proofs of the theorems in the supercuspidal case, and about the
restriction to $F=\qq_p$ in the results. As we have remarked in the abstract, and as should be clear from the statements of the theorems,
the results are based on reduction to characteristic $p$. This may seem somewhat artificial, but it may be motivated as follows.
As we have explained, the fact that $v_\mu$ is a test vector for the character $\chi$ boils down to showing that $v_\mu$
has nonzero projection on the $\chi$-eigenspace of $T$. In other words, one has to show that the product of the two idempotents
$e_\chi e_\mu$ is nonzero on $\rho(\nu)$, where $e_\mu=\frac{1}{q-1}\sum_{\sigma\in D^0} \mu^{-1}(\sigma)\sigma$, and
$e_\chi=\frac{1}{q^2-1}\sum_{\sigma\in T}\chi^{-1}(\sigma)\sigma$. The construction of Deligne-Luzstig (which is due to Drinfeld in this case)
shows how to realize the representation $\rho(\nu)$ (or rather, it's restriction to $SL_2(\ff_q)$) in the $\ell$-adic \'etale cohomology
of the smooth plane curve $C:X^qY-Y^qX =1$ over $\overline\ff_p$. 
Here the action of $SL_2(\ff_q)$ on $C$ is via the standard linear action on $(X, Y)$. Since $SL_2(\ff_q)$ acts on $C$, it also acts on the Jacobian
of $C$, and we find that there is a realization of $\rho(\nu)$ in the $\ell$-adic Tate module of $J$. Our problem is therefore to compute the image
of $e_\chi e_\mu$ in the $\rho(\nu)$-isotypic part of $T_\ell(J)\otimes\ol\qq_\ell$, and show that it is nonzero.

So far this idea is mere tautology -- computing $e_\mu e_\chi$ on the Jacobian is no different from trying to compute it directly, since the realization
via Deligne-Lusztig doesn't give any explicit model. However, the point is that one can think of $e_\mu e_\chi$ as a kind of endomorphism of $J$,
and to prove that an endomorphism $e$ is nonzero, it is enough to compute $e$ on the tangent space of $J$, or what amounts to the same, on the space of
holomorphic differentials on $C$. In other words, it is sufficient to check non-vanishing in the characteristic $p$ vector space given by 
de Rham cohomology rather than the characteristic zero \'etale cohomology. The de Rham cohomology of $J$ is a
representation of $SL_2(\ff_q)$ over $\overline\ff_p$, and one can hope that it is easier to compute with than the original representation in characteristic zero. 
Happily, this turns out to be just the case.

Let us assume henceforth that $q=p$ is prime.  
Since $C$ is a smooth plane curve of degree $p+1$, the classical method of adjoints shows that the sheaf of regular differentials is isomorphic to 
the twisting sheaf $O(p-2)$, and $H^0(C, \Omega^1)$ is the span of the space of monomials $X^iY^j$ with degree $d=i+j\leq p-2$. 
The action of $SL_2(\ff_p)$ is simply linear. Clearly the space of polynomials of a given degree is invariant under $SL_2(\ff_p)$, and 
we get a decomposition of $H^0(C, \Omega^1)$ in to the sum of $p-1$ representations, each of which is known to be irreducible.
These representations are essentially the set of possible Serre weights for $SL_2(\ff_p)$. 
The eigenvectors for the split torus $D^0$ are just the monomials $X^iY^j$. To find the eigenvectors of the nonsplit torus $T$, one has
 to diagonalize matrices of the form $\begin{pmatrix} a & b \\ b\xi& a\end{pmatrix}$, where $a, b\in\ff_p, a^2-\xi b^2\neq 0$ and $\xi$ is a nonsquare;
this can be accomplished by considering monomials in $U = X -\sqrt{\xi}Y, V= X + \sqrt{\xi}Y$, and the composite $e_\chi e_\mu$ can be computed
on any given polynomial simply by making the corresponding change of variables. 

This geometric approach was the one taken in an earlier version of this work. However, the details are somewhat involved, even though the underlying idea is very simple, since one has to pass somehow from $SL_2$ to $GL_2$, and keep track of the various representations in characteristic zero and characteristic $p$.
Furthermore, one has to deal with the fact that $e_\chi$ and $e_\mu$ will in general have non-rational coefficients, and are not genuine endomorphisms of the Jacobian.
Nevertheless, it is quite natural to compute with Deligne-Lusztig curves, given the realization of the discrete series of $SL_2(\ff_p)$
in the Drinfeld curve, and it is also natural in view of the appearance in our results of the weights in Serre's conjecture. This connection
 also arises in results of Herzig \cite{herzig}, who deals with the Serre weights which show up in the de Rham cohomology of Deligne-Lusztig
varieties for $GL_n$. In general, it seems like an interesting idea to compute the de Rham cohomology of a Deligne-Lusztig variety as a representation space for the group in question, and compare the results with those in characteristic zero. 

In this paper we use a  simpler and more elementary method: we simply reduce the characteristic zero representations modulo $p$, and compute $e_\chi e_\mu$
on the reduced representation. This turns out to be enough for our purposes, and avoids any complications coming from geometry. The reduction of a discrete
series representation of $GL_2(\ff_p)$ has generically $2$ components, and this gives rise to our Type 1/Type 2 dichotomy. 
However, implementing this method
reveals that it is significantly more complicated for residue field $\ff_q\neq \ff_p$ -- 
the reduction of a discrete series representation of $GL_2(\ff_q), q=p^d$ 
typically has 
$2^d$ components, parametrized by subsets of the $d$ embeddings of $\ff_q$ in $\ol\ff_p$. Thus we get correspondingly more types of characters: one has a notion
of Type 1/Type 2 at each embedding of $\ff_q$, and the statements of the results become correspondingly more complicated. We have not pursued this avenue,
in this paper, principally to keep the length of the article reasonable. However, it seems clear that the same method would work for general $q$, at least in principle.


I would like to thank Bill Casselman, Fred Diamond, Kimball Martin, and Rachel Ollivier for answering my questions, Dick Gross for his encouragement, and 
FLorian Herzig and Dipendra Prasad for some helpful comments on an early draft.

\subsection*{Notation:} Let $\oo$ denote the ring of integers in $F$, and let $P$ denote the maximal ideal, 
generated by the uniformizer $\varpi$. The ring of integers of $F$ is written as $\oo$. The residue field of $F$
is $k_F=\ff_q$. If considering a finite
extension $K/F$ we indicate the corresponding objects with a subscript $K$. 
For a ring $R$, let $G(R)$ denote the group $GL_2(R)$. 
For an integer $s\geq0$, let $G(s)$ denote the principal congruence
subgroup of $G(\oo)$, and let $B(s)$ denote the subgroup with
of matrices with lower left corner divisible by $P^s$. We let $T$ denote a maximal
torus of $G$, assumed non-split, and let $\Pi$ denote an irreducible infinite dimensional
representation of $G(F)$ which is ramified, in the sense that the conductor of $\Pi$ is a nontrivial
ideal of $\oo$. Equivalently, the space of invariants of $G(\oo)$ in $\Pi$ is trivial. If $K$ is any finite
extension of $F$ and $\chi$ is a character of $F^\times$, then the conductor of $\chi$ is the the
ideal $(\varpi_K^r)$, where $r$ is the smallest
non-negative integer such that $\chi$ is trivial on the set $1+ \varpi_K^r\oo_K$. We will sometimes
refer to $r$ itself as the conductor, and hope this will not ruffle too many feathers.

\section{$K$-types and the principal series}
\label{pseries}

In this section we use a key result due to Silberger Casselman, which states that any smooth irreducible admissible 
infinite dimensional representation of $GL_2(F)$ is the sum of a primitive piece of finite dimension (the ``type")
and another piece which depends only on the central character. Thus we will analyze test vectors
for principal series representations, and 
deduce the consequences for an arbitrary representation sharing the same central character. 

Consider an irreducible admissible smooth infinite dimensional representation of $GL_2(F)$
with central character $\omega$. Let $P^r$ denote the conductor of $\Pi$, so that 
$r$ is the smallest non-negative integer such that $\Pi$ contains a nonzero vector $v^\text{new}$
such that $g\cdot v^\text{new}=\omega(a) v^\text{new}$, for $g=\begin{pmatrix} a & b \\ c& d\end{pmatrix}
\in B(r)$.  For each positive integer $s\geq r$, we let 
$u_{P^s}(\omega)$ denote the representation definined by Casselman 
\cite{cass73}, page 312. Thus $u_{P^r}(\omega)=\ind^{G(\oo)}_{B(r)}(\omega)$, 
while for $s > r$ it is is defined
as the complement of $\sum_{ r\leq t < s}\ind_{B(t)}^{G(\oo)}(\omega)$ inside 
$\ind_{B(s)}^{G(\oo)}(\omega)$. 
Again, we view the character $\omega$ as a character of $B(s)$ via $\begin{pmatrix}
a & b \\c & d\end{pmatrix}\mapsto \omega(a)$, which makes sense as long as $s\geq r$. 
Each representation $u_{P^s}(\omega)$ irreducible, by Proposition 1(b) of \cite{cass73}. Each 
$u_{P^s}(\omega)$ contains a distinguished vector $u_s$ on which $B(s)$ acts via the character
$\omega$. The vector $u_s$ is unique up to scalars, and $u_{P^s}(\omega)$ is characterized as the unique
irreducible representation of $G(\oo)$ containing a $\omega$-eigenvector for $B(s)$ and which is trivial
on $G(s)$ but not on $G(s-1)$.  Then Casselman's theorem from {\cite{cass-restriction}} is the following:

\begin{thm} The restriction of $\Pi$ to $G(\oo)$ decomposes as $\Pi=V(\pi)\oplus
V(\omega)$, where $V(\pi)$ is finite dimensional, and $V(\omega)=\oplus_{s\geq r} u_{P^s}(\omega)$. 
The space $V(\pi)$ is the space of invariants under the principal congruence subgroup $G(r-1)$. 
\end{thm}

Consider a ramified principal series representation $\Pi$ of $GL_2(F)$ that has
minimal conductor amongst its twists. Thus $\Pi=\Pi(\nu, 1)$ for some character
$\nu:F^*\rightarrow \cc^*$. Let $P^r$ denote the conductor of $\nu$.
We assume that $r\geq1$. Then the conductor of $\Pi$ is the ideal $P^r$, which is
the same as the conductor of $\nu$. 

\begin{prp}
\label{cass}
 The restriction of $\Pi$ to $G(\oo)$ decomposes as $\oplus_{s\geq r}
u_{P^s}(\nu)$. 
\end{prp}

\begin{pf} It is clear that $\oplus_{s\geq r}u_{P^s}(\nu)$ occurs
inside $\Pi\vert_{G(\oo)}$. According to Theorem 1 in \cite{cass73}, the complement
is the space of vectors fixed by $G(r-1)$. However, this is zero, since the central
character of $\nu$ has conductor $P^r$. 
\end{pf}

\noindent Thus the first irreducible piece of $\Pi$ is the $(q+1)q^{r-1}$-dimensional 
representation $\ind^{G(\oo)}_{B(r)}(\nu)$, generated by the new vector.

Now consider a quadratic field extension $K/F$. We let $\oo_K$ denote the ring
of integers in $K$. We assume given an embedding of $\oo_K^*$ into $G(\oo)$. 
It is well-known that such embeddings exist, since $G(\oo)$ is the group of units in a
maximal order of $M_2(F)$. Thus, the restriction of $\Pi$ to $\oo_K^*$ factors through
the decomposition in the proposition above. Let $\oo_{K, s}$ denote the order
$\oo+ P^s\oo_K$, for a positive integer $s$.

\begin{lmm} We have $\oo_K\cap G(r)= \{u\in \oo_{K, r}|u\equiv 1\pmod {P^r}\}$. 
\end{lmm}
\begin{pf} Suppose $u\in \oo_K\cap G(r)$. Then $\frac{u-1}{\varpi^r}\in 
M_2(\oo)\cap K$, where $\varpi$ is a uniformizer for 
$P\subset\oo\subset \oo_K$. But $M_2(\oo)\cap K=\oo_K$, so we are done.
\end{pf}

Now we assume that $K/F$ is unramified. We want to restrict the decomposition given by
Casselman
above to the group $\oo_K^*\subset G(\oo)$, and for this we need some way to
label the characters of $O_K^*$ whose restriction to $O^*$ coincides
with the central character $\nu$ of $\Pi$. Consider any such character $\chi$. 
Define the conductor of $\chi$ to be the smallest positive integer $s$ such that $\chi$
is trivial on the set $\{u\in\oo_K^*\vert u\equiv 1\pmod{P^s}\}$. Clearly
we must have $s\geq r$. The order of $(\oo_K/P^s\oo_K)^*$ is 
$(q^2-1)q^{2(s-1)}$,
and the order of $(\oo/P^s)^*$ is $(q-1)q^{s-1}$, so there are $(q+1)q^{s-1}$
characters of conductor up to $s$. Note also that since $K/F$ is unramified, a character
of $K$ is determined uniquely by its restriction to $\oo_K^*$ and $F^*$.

\begin{prp} Let $s \geq r$. Then the restriction of $u_{P^s}(\nu)$ to $\oo_K^*$
is the sum of characters of level precisely $s$ that agree with $\nu$ on 
$F^*$. 
\end{prp}

\begin{pf} Start with the case of $s=r$. It follows from the lemma 
above that the restriction of $u_{P^r}(\nu)$ to $\oo_K^*$
is trivial on the set of elements congruent to $1$ modulo $P^r$. Equivalently,
this restriction is the sum of characters of conductor at most $P^r$, and
any character that shows up is equal to $\nu$ on the centre $F^*$. But the dimension
of $u_{P^r}(\nu)$ is $(q+1)q^{r-1}$, which is the number of possible characters. 
Since we know from a theorem of Tunnell that the restriction is multiplicity free,
it follows that the restriction of $u_{P^r}(\nu)$ to $\oo_K^*$ is precisely
the set of characters of $K^*$ of conductor $r$ that
agree with $\nu$ on the centre.

It remains to deal with the case of $s >r$. This is a simple induction,
using the fact that the dimension of $u_{P^s}$ is equal to the number of 
possible characters, together with the facts that the restriction is known to 
be multiplicity free, and that $u_{P^s}$ is trivial on $G(s)$. 
\end{pf}

The following theorem produces the requisite test vectors for minimal principal series representations
and unramified $K/F$. 

\begin{thm} Let $K/F$ be unramified. Let $\chi$
denote a character of $K^*$ such that $\chi$ agrees with $\nu$ on $F^*$
and such that $\chi$ has level $s> r$. Let $u_s$ denote a nonzero vector
on the unique line in $u_{P^s}(\nu)$ that is a $\nu$-eigenspace for $B(s)$. Then $u_s$ is a test
vector for $\chi$. If $s\leq r$, then the new vector $v_r$ is a test vector for $\chi$.
\end{thm}

\begin{pf} By virtue of the proposition above, we know where to look for test vectors for characters
of level $s$: we must look in the representation $u_{P^s}(\nu)$. Now, it is known that each
representation $u_{P^s}(\nu)$ contains a vector $u_s$ such that $g\cdot u_s=\nu(g)u_s$
for $g\in B(s)$. Furthermore, we have $G(\oo/P^s)=T(\oo/P^s)\cdot B(\oo/P^s)$, where 
$T(\oo/P^s)$ and $B(\oo/P^s)$ denote the images of $\oo_K^*$ and $B(s)$ in $G(\oo/P^s)$. (The latter claim 
can be checked by observing that $\oo_K^*\cap B(\oo)=\oo^*$, and by counting elements.) Thus, the translates of $v_s$ by elements of $\oo_K^*$ span the entire
space $u_{P^s}(\nu)$. Thus, if $\ell_\chi(v_s)=0$, we must have $\ell_\chi(\Pi(t)v_s)
=\chi(t)\ell_\chi(v_s)=0$ for all $t\in \oo_K^*$, and consequently, $\ell_\chi$ must
vanish identically on the entire space $u_{P^s}(\nu)$. But this is impossible, since
$\chi$ appears in $u_{P^s}(\nu)$.  The proof of the statement in the case $s\leq r$ is similar.
\end{pf}

\begin{rmk} In general, the vector $u_s$ is not a translate of the new vector, unless $r=s$. However, it is easy to check
that since $u_s$ is a test vector for $\chi$ of conductor $s$, then the translate $v_s = \begin{pmatrix} \varpi^{s-r} & \\ 0 & 1\end{pmatrix}$
is a test vector as well. We leave the details to the reader; the point is that $v_s$ has nonzero projection on to $u_s$, since it is not
contained in $\sum{t < s} u_{P^t}(\nu)$. 
\end{rmk}

\begin{txt}
Now we want to look at the case where $K/F$ is ramified. In this case, there are two
distinct maximal orders of $M_2(F)$ that contain $\oo_K$, and these two are swapped by the action
of a uniformizer of $\oo_K$. Equivalently, $\oo_K$ is contained in an Eichler order
$R$ of level $P$. We set $H=R^*$. 

Once again, we want to work out what sort of characters of $K$ show up in the irreducible
piece $u_{P^s}(\nu)$. 
Let us write $\varpi_K$ and $\varpi$ for the uniformizers of $\oo_K$ and $\oo$ respectively, and
consider a character $\chi$ of $K^*$ such that $\chi$
agrees with $\nu$ on $F^*$. Define the conductor of $\chi$ to be the smallest positive integer $t$
such that $\chi$ is trivial on the set $\{u\in\oo_K^*\vert u\equiv 1\pmod{\varpi_K^t\oo_K}\}$. Clearly,
we must have $t\geq 2r-1$, just by looking at the central character, which has conductor $r$ as a character
of $F^*$. 
\end{txt}

\begin{lmm} Let $\chi$ denote a character of $\oo_K^*$ such that $\chi$ agrees with $\nu$ on $\oo_K^*$. 
Then the conductor of $\nu$ is either equal to $2r-1$, or an even integer $2s$, where $s\geq r$.
\end{lmm}

\begin{pf} It is clear that the conductor of $\chi$ is at least $2r-1$, where $P^r$ is the conductor of $\nu$. 
Suppose that $\chi(u)=1$ for the set of elements $u\in\oo_K^*$ such that $u\equiv 1\pmod{\varpi_K^{2t+1}}$,
with $t\geq r$. Then since $\chi$ agrees with $\nu$ on $F^*$, $\chi$ is also trivial on the set of elements
$u'\in F^*$ such that $u'\equiv 1\pmod{\varpi^t}$.  Now if $u_1\in\oo_K^*$ is any element that is congruent
to $1$ modulo $p^{2t}$, we may write $u_1=1 + a\varpi_K^{2t} + b\varpi_K^{2t+1}$, with $a\in\oo$ and $b\in\oo_K$, since the residue
fields of $\oo_K$ and $\oo_F$ are equal. Now select $c\in\oo$ such that $c\varpi^t\equiv a\varpi_K^{2t}\pmod{\varpi_K^{2t+1}}$. Such
$c$ exists because $\varpi_K^{2t}$ and $\varpi^t$ are associate in $\oo_K$, and the residue
fields of $\oo_K$ and $\oo_F$ are equal. Then $\chi(1-c\varpi^t)=\nu(1-c\varpi^t)=1$, 
and $\chi(u_1) = \chi((1 - c\varpi^t)u_1)=1$, because $(1 - c\varpi^t)u_1\equiv 1\pmod{\varpi_K^{2t+1}}$. Thus the 
conductor of $\chi$ is at most $2t$, and we are done.
\end{pf}

 We can now enumerate the characters $\chi$ of $K^*$ which agree with $\nu$ on $F^*$. 
   
\begin{lmm} There are exactly $2q^{r-1}$ characters of $K^*$ 
of conductor $2r-1$ that agree with $\nu$ on $F^*$, and these restrict to $q^{r-1}$ distinct 
characters of $\oo_K^*$ of conductor $2r-1$ that agree with $\nu$ on $\oo^*$. 

If $s\geq r$, there are exactly $2(q^s-q^{s-1})$ characters of $K^*$ 
of conductor $2s$ that agree with $\nu$ on $F^*$, and these restrict to $q^s-q^{s-1}$ distinct 
characters of $\oo_K^*$ of conductor $s$ that agree with $\nu$ on $\oo^*$. \end{lmm}

\begin{pf} It suffices to count the number of distinct characters of $\oo_K^*$ of given conductor 
 which agree with $\nu$ on $\oo^*$, since each
such obviously extends in two different ways to $K^*$, according to the value
on the uniformizer $\varpi_K$, and the latter is determined up to sign by the restriction to $F^*\oo_K^*$. 

In the case of minimal conductor $2r-1$, we see that 
$(\oo_K^*/P^{2r-1}\oo_K)^*$ has cardinality $(q-1)q^{2r-2}$. If $r=1$, then 
the value of the character is determined by the values on elements of $\oo^*$, so it is unique.
If $r > 1$, then the value of the character is determined on the elements coming from $\oo^*$.
Observe that elements of $\oo$ have $\varpi_K$-adic expansions involving only even powers of $\varpi_K$, and
are of the form $a_0 + a_1\varpi_K^2 + \dots a_{r-1}\varpi_K^{2r-2}$ modulo $\varpi_K^{2r-1}$, which gives
$(q-1)q^{r-1}$ elements. Thus there are $q^{r-1}$ possible values for the character, if its value is fixed
on $\oo^*$. 

The remaining cases over even exponent $2s$ are calculated by a simple induction, starting with $r=s$. 
If $s\geq r$, then to calculate the number of characters of $\oo_K^*$ 
of conductor up to $2s$ that agree with $\nu$ on $\oo^*$,
we simply observe that $(\oo_K^*/P^s\oo_K)^*$ has cardinality $(q-1)q^{2s-1}$, whereas 
$(\oo/P^s\oo)^*$ has cardinality $(q-1)q^{s-1}$. 
\end{pf}

\begin{txt} We know abstractly that the representation $\Pi$ decomposes as the direct sum of all characters $\chi$
of $K^*$ which agree with $\nu$ on $F^*$, and that this restriction is multiplicity-free. 
Thus if we restrict further to a representation of $\oo_K^*$, we get each possible character with multiplicity 
two. If the maximal order $\oo_K$ is optimally embedded in the maximal order $M_2(\zz_p)$, then this decomposition
is compatible with the decomposition in (\ref{cass}), and our tasks is to locate the various eigenspaces for the full
torus $K^*$ in terms of (\ref{cass}). Each eigenspace for $\oo_K^*$ is a two dimensional subspace
inside  (\ref{cass}), although we shall see below that it is typically not contained in any single summand. 
The unformizer $\varpi_K$ 
acts nontrivially on every such eigenspace, and splits it in to two eigenspaces for distinct characters of the full
group $K^*$. However, $\varpi_K$ is not contained in $G(\oo)$, so we have to be a bit careful in determining
its action, with respect to the decomposition (\ref{cass}). 
\end{txt}

\begin{txt}  When $K/F$ is ramified the order $\oo_K$ is contained
in two distinct maximal orders of the matrix algebra. Thus $\oo_K$ is contained in an Eichler
order of level $p$, and we may assume that this order is the standard one
consisting of matrices whose lower left entry is divisible by $\varpi$. Denote this
order as $R(P)$. It will be handy in the sequel to have an explicit description of the embedding at hand, which
may be given as follows. Since $K/F$ is ramified and quadratic, it is tamely ramified if the residue
characteristic is odd. In this case, we may assume that $\varpi_K^2=\varpi$ is a uniformizer for
$\oo$, and the required embedding is given by
$$\varpi_K\mapsto \begin{pmatrix} 0 & 1 \\ \varpi & 0\end{pmatrix}.$$

In the case that the residue characteristic is two, we may still write $=F(\sqrt{\xi})$ for some $\xi$, where
$\xi$ is a non-square in $F$, and $\ord_{P}(\xi)= 0, 1$. If $\ord_P(\xi)=1$, then it is a uniformizer for $\oo$, and 
we may take the same embedding as above. If $\xi$ is a unit, then since we are assuming that $K/F$ is ramified,
this implies that $\ord_{P}(\xi-1)=1$ and $\varpi=\xi-1$ is a uniformizer for $\oo$. In this case, our embedding is given 
by 
$$\varpi_K \mapsto \begin{pmatrix} 1-\xi& -\xi\\\xi -1 & \xi -1\end{pmatrix}.$$

With these embeddings in hand, 
we will attempt to decompose the representation of $\ind^{G(\oo)}_{B(r)}(\nu)$ according to the characters
of $\oo_K^*$ that it contains. By definition, this representation is realized in the space
of translates of a new vector $v_r$, on which $B(r)$ acts via the character $\nu$, where 
we view $\nu$ as a character of $B(r)$ by evaluating $\begin{pmatrix}a & b \\c & d\end{pmatrix}
\mapsto\nu(a)$, as always. 
 \end{txt} 

\begin{lmm} The subgroup $1 + \varpi_K^{2r-1}\oo_K \subset B(r)$ acts trivially on $v_r$. The subspace of 
$\ind^{G(\oo)}_{B(r)}(\nu)$ generated by the vectors $t\cdot v_r$, for $t\in\oo_{K}^*$
has dimension $q^{r-1}$, and is the direct sum of distinct characters of $\oo_K^*$ of conductor $2r-1$.
\end{lmm}

\begin{pf} The first statement is a direct computation, using the embeddings given above, since
elements of $1 + \varpi_K^{2r-1}\oo_K$ are represented by matrices which have
have diagonal elements that are congruent to $1$ modulo $\varpi^r$, and $\nu$ is trivial on such 
elements. The second statement follows from the fact that $\oo_K\cap R(P^r)$ is the order 
$\oo_{K, r-1}$, and thus the vectors $t\cdot v_r$ for $t\in \oo_K^*/\oo_{K, r-1}^*$ 
are linearly independent in the \emph{irreducible} induced representation $\ind^{G(\oo)}_{B(r)}(\nu)$.
The fact that the characters of $\oo_K^*$ that appear are distinct follows from the linear
independence of the vectors $t\cdot v_r$, which implies that they realize a quotient of the regular representation
of $\oo_K^*/\oo_{K, r-1}^*$.
\end{pf}

Let $V_r$ denote the space spanned by the vectors  $t\cdot v_r$, for $t\in\oo_{K}^*$. Let
$w_r=\varpi_K\cdot v_r$, and let $W_r=\varpi_K\cdot V_r$. 

\begin{lmm}  The following assertions hold:
\begin{enumerate}
\item The space $W_r$ is spanned by the vectors $t\cdot w_r$.
\item The space $W_r$ is contained in the representation  $\ind^{G(\oo)}_{B(r)}(\nu)$. 
\item The spaces $V_r$ and $W_r$ generate a $2q^{r-1}$-dimensional space inside $\ind^{G(\oo)}_{B(r)}(\nu)$
in which each character of $\oo_K^*$ of conductor $2r-1$ which agrees with $\nu$ on $\oo_K^*$ appears
with multiplicity $2$. 
\end{enumerate}
\end{lmm}

\begin{pf} The first assertion is obvious, since $K^*$ is commutative. For the rest, 
we may argue as follows. According to Casselman, the space $\ind^{G(\oo)}_{B(r)}(\nu)$ is the space of invariants in $\Pi$
under the principal congruence subgroup $G(r)$. Thus it suffices to show that $w_r$ is fixed under
$G(r)$. Let $g\in G(r)$. Then direct computation with the matrix representing $\varpi_K$ shows that $\varpi_K^{-1} g \varpi_K$
is a matrix in $B(r+1)$, with diagonal entries that are 1 modulo $\varpi^r$. Since $w_r=\varpi_K\cdot v_r$, and $B(r+1)\subset
B(r)$ acts on $v_r$ via a character, the second statement of the lemma follows. 

As for the final statement, it follows from the fact that $\ind^{G(\oo)}_{B(r)}(\nu)$ is irreducible of dimension 
$(q+1)q^{r-1}$ that $\ind^{B(1)}_{B(r)}(\nu)$ is irreducible, of dimension $q^{r-1}$ as a representation of $B(1)$. 
Using the explicit
embeddings given above, once checks that $\varpi_K B(1)\varpi_K^{-1}=B(1)$, so $W_r$ is also a representation
space for $B(1)$, and is also irreducible. Thus either $V_r=W_r$ or $V_r\cap W_r=0$. We claim that the latter must
hold. To verify this, observe that 
the vector $w_r\in W_r$ is stable under $B_r$, and realizes the representation $\begin{pmatrix}a & b\\cp^r & d
\end{pmatrix}\mapsto \nu(d)$ (as opposed to $\nu(a)$, which picks out $v_r$). This can be checked directly, by using the specified embeddings. Thus, it suffices to prove that $\ind^{B(1)}_{B(r)}(\nu)$ contains no such vector. 
To do this, we have to find the dimension of the space of functions $f$
on $B(1)$ such that $f\left(\begin{pmatrix} a' & b'\\ 0 & d'\end{pmatrix}g \begin{pmatrix} a' & b'\\ 0 & d'\end{pmatrix}\right)= 
\nu(a')\nu(d)f(g)$. It is well-known that $B(r)\backslash B(1)/B(r)$ has $r$ elements, 
represented by the elements $\begin{pmatrix} 1 & 0 \\ \varpi^j & 1\end{pmatrix}$, for $1\leq j\leq r$, 
and using these representatives, we may compute as in \cite{cass73}, page 305, to see that no such function exists. 
Thus $V_r\cap W_r=0$, and we are done.
\end{pf}

Now we can find test vectors for the characters of $K^*$ with minimal level $2r-1$. We remind
the reader that we are assuming that $\oo_K$ is contained in the Eichler order $R(P)$,
and embedded as above.

\begin{prp} Let $\chi$ denote any character of conductor $2r-1$ on $\oo_K^*$ which agrees with
$\nu$ on $F^*$.  Let $v_r$ denote a nonzero vector such that $g\cdot v_r =\nu(g) v_r$, for $g\in B(r)=R(P^r)^*$ (i.e.
$v_r$ is a new vector). 
Then $v_r$ is a test vector for the unique $\chi$-linear form $\ell_\chi$ on $\Pi$.
\end{prp}

\begin{pf} Suppose that $\ell_\chi(v_r)=0$. Then $\ell_\chi(\Pi(t)v_r)=0$ as well, for all  $t\in K^*$. 
In particular, $\ell_\chi$ vanishes on the entire $2q^{r-1}$-dimensional
space spanned by $W_r$ and $V_r$, since they are generated by $K^*$-translates of $v_r$. 
But this is impossible, since the character $\chi$
occurs in this space.
\end{pf}

\begin{rmk} Note that we have not used any particular property of $\pi(\varpi_K)$, other than the fact
that $\varpi_K$ normalizes $B(1)$. One could of course use the fact that $\Pi$ is induced as a representation
of $G(F)$ to compute $\varpi_K\cdot v_r$ explicitly, but we prefer not to do this, as we want an argument
that depends only on the components of the restriction of $\Pi$ to $G(\oo)$.
\end{rmk}

It is now easy to see how to get test vectors for the characters of next level up. The representation
$\ind^{G(\oo)}_{B(r)}(\nu)$ has dimension $(q+1)q^{r-1}$, and there is a $2q^{r-1}$-dimensional subspace on which 
$\oo_K^*$ acts via characters of conductor $2r-1$. Let $V_{r+1}$ denote the complementary $\oo_K^*$-invariant subspace, so that $V_{r+1}$ has dimension $q^{r}-q^{r-1}$. 

\begin{lmm} The representation of $\oo_K^*$ on $V_{r+1}$ is given as the direct sum of characters $\chi$, where
$\chi$ runs over the characters of conductor $2r$ which agree with $\nu$ on 
$F^*$. In particular, each such character occurs with multiplicity $1$.
\end{lmm}

\begin{pf} Consider the tree of $PGL_2$, whose vertices correspond to maximal compact subgroups of $G(F)$.The maximal
subgroup $G(\oo)$ defines a vertex $x$ in the tree. This vertex has $p+1$ neighbours, precisely one of which contains
$O_K^*$, namely, the other endpoint of the edge corresponding to the standard Eichler order $R(P)$. 
Consider a vertex $y$  at distance $r$ from $x$, with 
corresponding maximal compact $K_y$. Then $\oo_K\cap K_y =\oo_{K,r}^* $ or $\oo_K\cap K_y=\oo_{K,{r-1}}^*$, 
depending on whether or not $K_y$ contains the units of the standard Eichler order
$R(P)$. Equivalently, $\oo_K\cap K_y=\oo_{K,{r-1}}^*$ if and only if the path from $y$ to $x$ contains the edge
corresponding to $R(P)$. Thus, there are $q^{r-1}$ vertices $y$ which satisfy $\oo_K\cap K_y=\oo_{K,{r-1}}^*$,
and $q^r$ that satisfy $\oo_K\cap K_y=\oo_{K,r}^*$. Let $R_y$ denote the Eichler order of discriminant $P^r$ whose 
unit group is given
by the intersection of $K_x$ and $K_y$. The orders $R_y$ are permuted transitively by the conjugation action of $G(\oo)$,
since $G(\oo)$ acts transitively on the set of maximal orders at distance $r$. The line spanned by the 
vector $v_r$ appearing previously is fixed by the units of the standard order
$R(P^r)$ whose lower left entries are divisible by $P^r$, which corresponds to $R_{y_0}$
for some vertex $y_0$ at distance $r$. It follows that the unit group of each order $R_y$ fixes a line spanned by a vector
 $v_y\in \ind^{G(\oo)}_{B(r)}(\nu)$. Since $\ind^{G(\oo)}_{B(r)}(\nu)$ has dimension $(q+1)q^{r-1}$, the span of the vectors
 $v_y$ is stable under $G(\oo)$, it follows that the vectors
 $v_y$ are linearly independent and span $\ind^{G(\oo)}_{B(r)}(\nu)$.  
The space $V_r$ mentioned previously is the span of the vectors $v_y$ for the $q^{r-1}$ orders $R_y$ such that 
$\oo_K\cap K_y=\oo_{K,r-1}^*$. 

Let  $y$ be such that $\oo_K\cap K_y=\oo_{K,r}$. Let $v_y$ be a nonzero vector spanning a line fixed by $R_y^*$. 
Then the $\oo_K^*$-orbit of $v_y$ is spanned by the vectors $v_{y'}$, corresponding to the lines fixed by $R_{y'}^*$,
where $R_{y'}$ runs over the $\oo_K^*$-orbit of $R_y$. But it is clear that this orbit has cardinality given
by $\oo_{K}^*/\oo_{K, r}^* = (q-1)q^{2r-1}/(q-1)q^{r-1}=q^r$.  It follows therefore that the space of
$\oo_K$-translates of $v_{y}$ generate the representation of $\oo_K^*$ induced from a character of $\oo_{K, r}^*=
\oo_K\cap R_y^*$.
Thus we get all characters of $\oo_K^*$ of conductor up to $2r$ which agree with $\nu$ on the centre, 
each with multiplicity $1$. The space $W_r$ accounts for $q^{r-1}$ characters of conductor $2r-1$, and the remaining
$q^r-q^{r-1}$ characters give the space $V_{r+1}$, containing the ones of conductor $2r$, each with multiplicity one.
\end{pf}

We may now repeat the same argument as was used in the minimal case, finding a disjoint space 
$W_{r+1}\subset \ind^{G(\oo)}_{B(r+1)}$ consisting of vectors
spanning lines fixed by orders at distance $r$ from the other maximal order which contains $\oo_K^*$, 
and at distance $r+1$ from $K_x$. The two
spaces are interchanged by the action of the uniformizer $\varpi_K$ on $V_{r+1}$ to obtain a 
disjoint space which  realizes the same set of characters , and the same argument
then shows that the vector $v_y$, where $y$ is any vertex at distance $r$ from $x$ such that $R_x$ optimally contains
$\oo_{K, r}$ is a test vector for the characters of conductor $2r$. Repeating this process, we get the following result:

\begin{thm} Let $s \geq r$. Let $\chi$ denote any character of $K^*$ of conductor $2s$ which agrees with 
$\nu$ on $F^*$.  Then if $S$ is an Eichler order of discriminant
$P^s$ which optimally contains the order of conductor $P^s$, then any nonzero vector on the line fixed by $S$ is a test vector for $\chi$.
\end{thm}

\subsection{Consequences: test vectors for sufficiently ramified characters $\chi$}

We now take up the case of a general smooth irreducible admissible infinite dimensional representation $\Pi$ with central character $\omega$.
Let $P^r$ denote the conductor of $\omega$. Then we have $\Pi=V(\pi) \oplus_{s\geq r} u_{P^s}(\omega)$,
by Casselman's theorem again. Each $u_{P^s}(\omega)$ contains a nonzero vector $v_s$, where $V_{r_1+r_2}=v^{\text{new}}$ is the new vector, and $v_s$
is the translate of $v^{\text{new}}$ by the matrix $\begin{pmatrix} p^t & 0 \\ 0 & 1\end{pmatrix}$ for $t=s-(r_1+r_2)$. 
The analysis of $u_{P^s}$ made above shows immediately that $v_s$ is a test vector for any $\chi$ occurring in $u_{P^s}$. Thus we get the
following key theorem, which was proven by a more complicated method in \cite{fmp}:

\begin{thm} 
\label{fmpthm}
Let $\Pi$ be as above, and let $\chi$ be any character of $T$ that occurs in $u_{P^s}(\omega)$ for  $s\geq r$. Then $v_s$ is a test vector
for $\chi$. If $T$ is unramified, $\chi$ occurs in $u_{P^s}(\omega)$ if and only if $\chi$ has conductor $P^s$. Thus if $\chi$ is any character
of $T$ which agrees with $\omega$ on $F^\times$, and $\chi$ has conductor $P^s$ with $s\geq r$, where $P^r$ is the conductor of $\omega$,
then $v_s$ is a test vector for $\chi$.
\end{thm} 

In view of this result, the test vector problem will be completely solved if we can find test vectors for the characters $\chi$ 
occurring in the primitive part $V(\pi)$. We start with the case of a general principal series 
representation $\Pi=\Pi(\nu_1, \nu_2)$, where both $\nu_i$ are nontrivial. The central character
is $\omega=\nu_1\nu_2$. Let $\chi$ denote a character of $K^\times$ where $K/\qq_p$ is a quadratic field extension, possibly ramified. Assume that
$\chi$ and $\omega$ agree on $F^\times$. Let $r$ denote the smallest positive integer such that $\omega$ is trivial on $1+\varpi^r\oo_F$. Let $r_1$ and
$r_2$ denote the corresponding quantities for $\nu_1$ and $\nu_2$, and that $r_1\geq r_2$.
Consider $\Pi'=\Pi\otimes\nu_2^{-1}$. Then $\Pi'=
\Pi(\nu_1/\nu_2, 1)$ is also irreducible, then the results proved above apply to $\Pi'$. Let $P^{r'}$ denote
the conductor of $\omega'=\nu_1/\nu_2$. Then for $s\geq r'$, we have a vector $v'_s\in u_{P^s}(\omega')\subset \Pi'$, 
and it is clear that $w_s=v'_s\otimes\nu_2 \in u_{P^s}(\omega')\otimes\nu_2\subset\Pi$
is a test vector for any character $\chi$ such that $\chi\otimes\nu_2^{-1}$ occurs in $u_{P^s}(\Omega')$. This covers \emph{all} possible
$\chi$, since $\Pi'=\oplus u_{P^s}(\omega')$, according to the analysis carried out above. In particular, we obtain test vectors for characters
$\chi$ occurring in the primitive part $V(\Pi)$, which leads to the following general theorem:

\begin{thm}
\label{twisted} Suppose that $\Pi=\Pi(\nu_1, \nu_2)$ is as above. Let $\chi$ denote any character of $T$ which agrees with $\omega=\nu_1\nu_2$ on 
$F^\times$. Then, a test vector for $\chi$ is given by the vector $w_s=v'_s\otimes\nu_2$, where $v'_s\in u_{P^s}(\omega')$ is the test vector at
level $s$ for $\Pi'=\Pi(\nu_1/\nu_2, 1)=\oplus u_{P^s}(\omega')$, and $\omega=\nu_1/\nu_2$.
\end{thm}

\begin{rmk} Observe that if $s \geq r$, then we get test vectors in the case covered by Theorem \ref{fmpthm}, which are distinct from
the test vectors provided by that result, as can be seem by considering the action of the diagonal torus. 
\end{rmk}

The analogue of Theorem \ref{twisted} for special representations is holds as well.  In the case of the `minimal' special representation of trivial
central character, the space
$V(\Pi)$ is zero and the result of Theorem \ref{fmpthm} is comprehensive.

Thus, we have almost completely solved the test vector problem for special representations and the irreducible principle series, and a non-split 
torus $T$. The only question remaining is the following: if $\Pi$ is a non-minimal principal series, and $\chi$ is a character of $T$
which occurs in the primitive part $V(\Pi)$, then can we find a test vector for $\chi$ of the form $\begin{pmatrix} p^t & 0\\ 0 & 1\end{pmatrix}\cdot
v^\text{new}$, for some $t$ (necessarily negative)? 

As for the supercuspidal representations, Theorem \ref{fmpthm} gives the result for characters of conductor divisible by $P^r$, where $r$
is the conductor of $\Pi$. But there is no immediate way to deal with the characters of conductor $P^s$ with $s\leq r$. In particular,
even the \emph{trivial} character poses problems. (We remind the reader that the Gross-Prasad test vector for the trivial character is not 
fixed by the units in any Eichler order.)  One can guess
that test vectors of the form $\begin{pmatrix} p^t & 0\\ 0 & 1\end{pmatrix}\cdot v^text{new}$ exist, but, as we shall see in the next section,
 the answer is surprisingly delicate, even in the case of depth zero supercuspidals of conductor $p^2$ of $GL_2(\qq_p)$.

\section{The depth zero supercuspidal representations of $GL_2(F)$.} 
\label{supercusp}

Our goal in this section is to completely solve the test vector problem for depth zero supercuspidal representations of $GL_2(F)$ when $F=\qq_p$. 
and the torus $T$ corresponding to an unramified quadratic extension $K/\qq_p$. In this entire section, we assume that $p$ is odd.
  
We start by defining the groups and representations of interest. Let $G=GL_2(\ff_q)$, where $\ff_q$ is the residue field of $\oo=\oo_F$. Let $\nu$ denote
a character $\nu:\ff_{q^2}^*\rightarrow \cc^*$ which does \emph{not} factor through the norm $\ff_{q^2}\rightarrow\ff_q$, and let $\rho=\rho(\nu)$
denote the cuspidal representation of $G=GL_2(\ff_q)$ associated to $\nu$ by Piatetski-Shapiro \cite{ps-book}. Then there exists a representation
$\Pi=\Pi(\nu)$ of $GL_2(F)$ such that the type of $\Pi$ is the inflation of  $\rho$ to $GL_2(\oo)$. The representation $\Pi$ is called the depth 
zero supercuspidal representation associated to $\nu$. Our goal is to find a test vector for characters $\chi$ of the unramified torus $K/F$ and the representation 
$\Pi=\Pi(\nu)$ when $\chi$ occurs in the type, in the following sense. 

Fix $\Pi=\Pi(\nu)$, with central character $\omega$. 
As we have observed already, $T= \varpi^\zz \times U$ where $U\subset GL_2(\oo)$ is compact and (since $T$ is unramified) $\varpi$ is in $F^*$. If $\chi$
is a character of $T$ which agrees with $\omega$ on $F^*$, then $\chi(\varpi)$ is fixed, and $\chi$ is determined by its restriction to $U$. On the other hand,
Casselman's theorem tells us that
the restriction of $\Pi(\nu)$ to $GL_2(\oo)$ is the sum of the inflated representation $\rho(\nu)$ (the type) and the imprimitive part, so we get a corresponding
decomposition of $\Pi(\nu)$ as an infinite direct sum of characters of $U$. We say that $\chi$ occurs in the type if there is a line  in $\rho(\nu)$  on which $T$ acts by the character $\chi$. Our goal is to find a test vector for such $\chi$, since the characters occurring in the imprimitive part may be dealt with by the arguments in the first section of this paper. The point is that any $\chi$ which appears in the type is trivial on the principal units $u_1$ in $T$, since the type is inflated from a representation of 
$GL_2(\ff_q)$, which is trivial on $I + \varpi M_2(\oo)$. In other words, such a character $\chi$ may be identified with a character of $T(\ff_q)=U/U_1\equiv \ff_{q^2}\subset 
GL_2(\ff_q)$, and finding a test vector for such $\chi$ is reduced to a problem in finite group theory for the representation $\rho(\nu)$ of the group $GL_2(\ff_q)$. . 

As for the candidate test vectors, we are looking for vectors related to local newforms. By Atkin-Lehner theory, the local newform $v^\text{new}\in \Pi(\nu)$
is a nonzero vector on the unique line in $\Pi(\nu)$ where the group $B(2)$ acts via the character 
$\begin{pmatrix} a & b \\ p^2c & d\end{pmatrix}\mapsto \omega(a)$. 
As explained in Section 1, this means that $v^\text{new}\in \Pi(\nu)$ lies in the first imprimitive 
piece of $\Pi(\nu)$. On the other hand, consider $w^{new}=
\begin{pmatrix} \varpi & 0 \\ 0 & 1\end{pmatrix}\cdot v^\text{new}$, where the action of the matrix on $v^\text{new}$ 
is given by the representation $\Pi(\nu)$. Then
$w^\text{new}$ is invariant under $I + \varpi M_2(\oo)$, and therefore lies in the type. Furthermore, $w^\text{new}$ is an eigenvector for the matrices
$\begin{pmatrix} a & b\varpi \\ c\varpi & d\end{pmatrix}$, which is to say, they are eigenvectors for the diagonal \emph{split} 
diagonal torus $D$ of $GL_2(\ff_q)$. 
We want to know whether or not the unique $\chi$-equivariant functional on $\Pi(\nu)$ is nonzero on $w^\text{new}$, where $\chi$ is \a character of $T$ such
that $\chi=\omega$ on $F^*$ and 
such that the restriction of $\chi$ to $U$ occurs in the representation $\rho(\nu)$ of the group $GL_2(\ff_q)$. 
As remarked above, the fixed central character
allows us to identify characters of $T$ which agree with $\omega$ on $F^*$ with characters of $T(\ff_q)\subset GL_2(\ff_q)$.

Thus, we are led to
the following problem in the representation theory of the finite group $GL_2(\ff_q)$: suppose $\rho(\nu)$ is a cuspidal  representation of $GL_2(\ff_q)$. Suppose that
$\chi$ is a character of the \emph{nonsplit} torus $T(\ff_q)$ which occurs in the restriction of $\rho(\nu)$ to $T(\ff_q)$, and $\mu$ a character of the \emph{split}
torus $D$ which occurs in the restriction of $\rho(\nu)$ to $D$.  Let $v_\chi$ and $v_\mu$ denote corresponding eigenvectors. 
We will see below that each of these characters $\chi, \mu$ occurs with multiplicity one, and it follows
that the unique $\chi$-equivariant functional on $\rho(\nu)$ factors through the orthogonal projection to the line spanned by $v_\chi$.  Thus, to determine whether or not $v_\mu$ is a test vector for $\chi$ is equivalent to determining
whether or not the orthogonal projection of 
$v_\mu$ on to $v_\chi$ is nonzero. Here orthogonality is taken with respect to a $GL_2(\ff_q)$-invariant inner product.

\subsection{Discrete series representations of $GL_2(\ff_q)$ and test vectors modulo $p$} 

 In this subsection (and in this subsection only) we will use slightly different notation
from the rest of this paper, since we are working solely with representations of $GL_2(\ff_q)$. Our goal to solve the problem described in
the paragraph above.

\begin{txt}
Thus Let $\xi\in \ff_q^*$ denote a non-square element, and write $T$ for the subgroup of $GL_2(\ff_q)$ given by matrices of the
form $\begin{pmatrix} a & b \\b\xi & a\end{pmatrix}$, where $a, b\in \ff_p, a\neq 0$, so that
$T\cong \ff_{q^2}^*$. Let $D\subset S$ denote the subgroup $\begin{pmatrix} a & 0\\0& b\end{pmatrix}$ with $ a, b \in\ff_q^*$. 
Let $\omega:\ff_q^*\rightarrow\cc^*$ denote the central character of $\rho$. 
Thus $\omega$ is the restriction of $\nu$ to $\ff_q^*\subset\ff_{q^2}^*$. 

Next we consider the restriction of $\rho$ to the subgroups $D$ and $T$. 
The restriction of $\rho$ to $D$ consists of the direct sum of the $q-1$ characters 
$\mu:D\rightarrow{\cc}^*$ such that $\mu$ agrees with $\nu$ on the centre, each with multiplicity $1$. 
The characters $\mu$ can be distinguished by their restrictions
to the group $D^0=\begin{pmatrix} a & 0\\0& 1\end{pmatrix}, a\in \ff_p^*$. We write $v_\text{triv}$ 
and $v_\text{quad}$ for the corresponding eigenvectors for the action of $D^0$.

On the  other hand, the restriction of $\rho$ to $T$ is the sum
of $q-1$ characters  $\chi:T\rightarrow {\cc}^*$ which agree with $\nu$ on the centre; there are $q+1$ such characters, and each character other than 
$\nu$ and $\nu^q$ shows up with multiplicity one. Note that $\nu^q\neq \nu$, since $\nu$ does not factor through
the norm. In this situation, we shall say that a vector $v\in\rho(\nu)$ is a test vector for $\chi$ if 
the projection of $v$ to the $\chi$-eigenspace of $T$ is nonzero.  
\end{txt}

To state the general results, 
it is important to have a way to label the various characters of $\ff_{q^2}^*$. This is quite complicated
for a general finite field, so from now on, we assume that $q=p$ (and as always that $p$ is odd).  

\begin{txt}  We recall the labelling of characters of $\ff_{p^2}^*$ that was given in the introduction. 
If $x\in\ff_{p^2}$, 
we identify $x$ with its
image under the Teichm\"uller lift to the ring $W=W(\ff_{p^2})$ of Witt vectors for 
$\ff_{p^2}$. Let $\ol{x}=x^p$. Then $x\mapsto \ol{x}$ is also a
$W$-valued character of $\ff_{p^2}^*$. If $\chi$ is an arbitrary character of $\ff_{p^2}^*$ with values in $W$, we may write $\chi(x)=x^a\ol{x}^b$,
for positive integers $a, b$ satisfying $0\leq a, b \leq p-1$, since $\chi$ has order dividing $p^2-1$. We call this character $\chi(a, b)$.
 Note that we can exclude $a=b=p-1$, since
$a+pb$ may be assumed to lie in the range $0\leq a + pb < p^2-1$. 
We set $k_\chi= a + pb$, so that $\chi(x)=x^{k_\chi}$. 
We say that $\chi$ is even or odd according to the sign of $\chi(-1)$; equivalently,
$\chi$ is even if and only if $a, b$ have the same parity. The character $\chi$ factors
through the norm if and only if $a=b$, or -- equivalently -- if $k_\chi$ is divisible by $p+1$.  
Note that switching $\nu$ and $\nu^p$ interchanges $a$ and $b$,
so we may assume that $a > b$, if $\nu$ does not factor through the norm. 
\end{txt}

\begin{txt} 
An important point for us will be to determine the various characters $\nu$ that coincide on $\ff_p^*$:
 characters $\chi, \nu$ which agree on $\ff_p^*$ are those with weights $k_\chi\equiv k_\nu\pmod{p-1}$. If $\nu$ is a fixed
character with $\nu\neq \nu^p$, and $\chi$ is a character which agrees with $\nu$ on $\ff_p^*$, then we say that $\chi$ is of Type 1
if $k_\chi$ is between $k_\nu$ and $k_{\nu^p}$; else we say $\chi$ is of Type 2. In any case, we define $t_\chi=\frac{k_\chi - k_\nu}{p-1}$
Note that $\frac{k_{\nu^p} - k_\nu}{p-1}$ is even, so that $t_\nu$ and $t_{\nu^p}$ have the same parity. 
\end{txt}

With this enumeration of characters in hand, the result in the even case may be stated as follows. 
Consider the representation $\rho(\nu)$, where $\nu=\nu(a, b)$ is even. Let $\omega$ denote the central character
of $\rho(\nu)$, so $\omega(x)=x^{a+b}$. Define characters $\mu_1, \mu_2$ according to the following prescription: 
$\mu_1(x)=x^{(a+b)/2}$, and 
$\mu_2(x)=x^{(a+b + p-1)/2}$. Further define
 $\mu_3(x)=x^{(a+b)/2-1}$ and $\mu_4(x)=x^{(a+b+p-1)/2-1}$. 
Let $v_{\mu_i}$ denote an eigenvector for $D^0$ with eigenvalue $\mu_i$. Note
that $\mu_1, \mu_2$ are the two square roots of the central character $\omega$. 
\begin{thm} 
\label{evencharcase}
Consider the representation $\rho(\nu)$ for $\nu(a, b)$, with notation as above. Assume $a, b$ are labelled so that $a > b$, and that $a-b$ is even. 
Let $\chi:T\rightarrow W^*$ denote a character such that $\chi=\nu$ on $\ff_{p}^*$, and $\chi\neq \nu, \nu^q$.Then the following 
statements hold:
\begin{enumerate}
\item Suppose $\chi$ is of Type 1 and $t_\chi$ is odd. Then the vector $v_{\mu_1}$ is a test vector for $\chi$. 
\item Suppose $\chi$ is of Type 2 and $t_\chi$ is odd. Then the vector $v_{\mu_2}$ is a test vector for $\chi$.
\item Suppose $\chi$ is of Type 1 and $t_\chi$ is even. Then the vector $v_{\mu_3}$ is a test vector for $\chi$. 
\item Suppose $\chi$ is of Type 2 and $t_\chi$ is even. Then the vector $v_{\mu_4}$ is a test vector for $\chi$.
\end{enumerate}
\end{thm}

Now consider the odd case, where $\nu=\nu(a, b)$ with $a  - b = 2t+1$.  Define characters $\mu_1,\mu_2$ of $\ff_p^*$ by $\mu_1(x) = x^{(a+b-1)/2}$ and
$\mu_2(x)=x^{(a+b+p)/2}$, and let $v_{\mu_i}$ denote eigenvectors for $D^0$ with eigenvalue $\mu_i$. Let
$\omega(x)=x^{a+b}$, for $x\in\ff_p^*$. If $\chi$ is a character of $\ff_{p^2}^*$ which agrees with $\omega$ 
on $\ff_p^*$, we define the notion of Type 1 and Type 2 exactly as before.

\begin{thm} 
\label{oddcharcase}
Consider the representation $\rho(\nu)$ for $\nu(a, b)$ odd, with notation as above. Assume $a, b$ are labelled so that $a > b$ and $a-b$ is odd. 
Let $\chi:T\rightarrow W^*$ denote a character such that $\chi=\omega$ on $\ff_{p}^*$, and $\chi\neq \nu, \nu^p$.
Suppose $a-b > 1$. Then the following 
statements hold:
\begin{itemize}
\item  Suppose $\chi$ is of Type 1. Then the vector $v_{\mu_1}$ is a test vector for $\chi$. 
\item Suppose $\chi$ is of Type 2. Then the
 vector $v_{\mu_2}$ is a test vector for $\chi$.
\end{itemize}

If $a-b=1$, then the vector $v_{\mu_2}$ is a test vector for $\chi$, for all $\chi$. 
\end{thm}

Once correctly formulated, these theorems are surprisingly easy to prove. 
The key point is that a cuspidal representation of $GL_2(\ff_p)$ 
is typically reducible modulo a prime above $p$, with Jordan-Holder filtration of length 2. 
The characters in question arise in one component
or the other, hence the Type 1 and Type 2 of the theorems. There is one exception case of odd central 
character, where the representation is irreducible
modulo $p$, leading to the final case of Theorem \ref{oddcharcase}.

\begin{txt} We start with the proof in the even case as stated in Theorem \ref{evencharcase}. 
Let $\nu$ be a character as in the theorem, such that $\rho(\nu)$ has even
central character. Let $L/F$ denote a finite extension of $F$ such that $\rho(\nu)$
may be realized via matrices with entries in $L$. There exists an $\oo_L$-lattice $B$ stable
under $\rho(\nu)$, where $\oo_L$ is the ring of integers in $L$, and the reduction 
of $B$ modulo the maximal ideal $\frak{m}_L$ of $\oo_L$ gives a $q-1$-dimensional 
representation $\ol\rho(\nu)$ of $GL_2(\ff_p)$ with
entries in $k_L$, where $k_L$ is the residue field of $\oo_L$. It follows from
\cite{prasad-modular}, Lemma 4.2, that this representation has 
Jordan-H\"older length 2, with composition factors $V_1=\det^{b+1}\operatorname{Sym}^{a-b-2}(k_L^2)$
and $V_2=\det^{a}\operatorname{Sym}^{p-1-(a-b)}(k_L^2)$.

Consider the representations of the two tori $T$ and $D^0$ on the $\ol\ff_p$-representations $\ol\rho(\nu)$ 
and  $V_i$. Each
$V_i$, as well as $\ol\rho(\nu)$, 
is semisimple for the actions of $D^0$ and $T$ separately, since these groups have order prime to $p$. 
Since each character in the decomposition of $\rho(\nu)$ with respect to $D$ or $T$ occurs with multiplicity
one, and since the groups have order prime to $p$, it follows that the same is true for the decomposition
of the reduced representation $\ol\rho(\nu)$. Here we are identifying $\ol\ff_p$-valued characters with their
Teichm\"uller lifts.
It follows that each of the $V_i$ may be decomposed as the sum of certain characters of either $T$ or $D$,
each with multiplicity one,
and that we get partitions of the characters of $D$ and $T$ occurring in $\rho(\nu)$, according to
whether they occur in $V_1$ or $V_2$. 

If $\chi$ is a character of $T$ and $\mu$ a character of $D$, we say that $\chi$ and $\mu$ are compatible
if they occur in the same component $V_i$ or not. If $v_\mu$ is a vector in the lattice 
$B$ which is an eigenvector for $D^0$ with eigenvector $\mu$, and $\chi$ character of $T$, 
then we shall say $v_\mu$ is a test vector for $\chi$ if $e_\chi(v_\mu)$ is nonzero modulo
modulo $\frak{m}_LB$, namely, nonzero in the characteristic $p$ representation obtained by
reduction of $B$. Clearly, a \emph{neccessary} condition for this to occur is that $\chi$ and
$\mu$ be compatible in the above sense. 

Our task is therefore to determine the decomposition of the representations $V_i$ to 
determine compatibility, and then to produce further conditions that exhibit test vectors
(modulo $p$) in each component. The characters listed in cases (1) and (3) of the 
Theorem \ref{evencharcase} are the characters
$\chi$ which occur in $V_1$, while the characters $\chi$ of cases (2) and (4) are the ones
occuring in $V_2$. 

Consider first $V_1=\det^{b+1}\operatorname{Sym}^{a-b-2}(k_L^2)$, which we realize as
the space of homogeneous 
polynomials in 2 variables $X, Y$ of degree $d= a-b-2$. Observe here that since $a$ and
$b$ are of the same parity, $a-b$ is even and positive. The group $GL_(\ff_p)$ acts linearly on $X, Y$. 
Each monomial $X^iY^{d-i}$ is an eigenvector for $D^0$, with eigenvalues given by the character $x\mapsto 
X^{i}$, for $b+1\leq i\leq a-1$.  Thus $e_\mu$ simply picks out the coefficient of $X^iY^{d-i}$, for 
suitable $i$ depending on $\mu$. The character $\mu_1$ in the theorem corresponds 
to the monomial 
$X^{d/2}Y^{d/2}$, while $\mu_3$ corresponds to $X^{d/2-1}Y^{d/2+1}$.

 To start with we compute 
 $e_\chi$ on the monomial $X^{d/2}Y^{d/2}$. Observe that $T$ is a nonsplit torus
 in $GL_2(\ff_p)$, and to find the eigenvectors for $T$, we must diagonalize the corresponding matrices
 $\begin{pmatrix} a & b \\ b\xi & a\end{pmatrix}$. Evidently, the eigenvectors are the monomials
 in variables $U=X-\sqrt{\xi} Y$ and $V=X+\sqrt{\xi}Y$, where $\sqrt{\xi}$ is some element of trace zero
 whose square is the nonsquare element $\xi$. Obviously, we have $X=\frac{1}{2}(U+V)$ and 
 $Y=\frac{-1}{2\sqrt{\xi}}(U-V)$, so that $X^{d/2}Y^{d/2}=
(\frac{-1}{2\xi})^{d/2}(U^2-V^2)^{d/2}=(\frac{-1}{2\xi})^{d/2}\sum_i \binom{d/2}{i} U^{2i} V^{d-2i}$. Each monomial
 is distinct, and appears with coefficient which is nonzero modulo $p$. It follows that 
 $e_\chi(X^{d/2}Y^{d/2})$ is nonzero for each character $\chi$ corresponding to a monomial in $U, V$
 which appears in the expression for $X^{d/2}Y^{d/2}$. The action of $x\in T$ on $U^{2i}V^{d-2i}$ is given
 by $x \mapsto x^{(p+1)(b+1)+ 2i + p(d-2i)}$, and after a brief calculation, one finds that the exponent is
$ (p+1)(b+1)+ 2i + p(d-2i) = b + ap -(2i+1)(p-1)$. Now the value of this exponent is $b+ap -(p-1)=k_{\nu^p}-(p-1)$ 
when $i=0$,  and $b + ap -(p-1)(a-b-1) = k_\nu + (p-1)$ when $i=d/2$, so we find that the characters
$\chi$ that show up are the ones with $k_\chi$ between $k_\nu$ and $k_{\nu^p}$, and of the form $k_\nu + (2j+1)(p-1)$,
as asserted. 
 
Moving on to Case 3, which is the other instance of the component $V_1$, consider this time
the monomial $X^{d/2-1}Y^{d/2+1}$. Making the change of variable $V=X+\sqrt{\xi}Y$, $U=X-\sqrt{\xi}Y$
as before, we find that  $X^{d/2-1}Y^{d/2+1}= C\cdot (U^2-\xi V^2)^{d/2-1}(U-V)^2=
C\cdot (U^2-\xi V^2)^{d/2-1}(U^2-2UV + V^2)$. Since $(U^2-\xi V^2)^{d/2-1}$ is a sum
of $d/2$ distinct monomials (the binomial coefficients being prime to $p$) with $U, V$
occurring to even powers, the expression $C\cdot (U^2-\xi V^2)^{d/2-1}(U^2-2UV + V^2)$ contains
$d/2$ distinct monomials where the exponent of both $U$ and $V$ is odd. The corresponding
characters $\chi$ are the ones listed in case 4 of the theorem, by calculating the corresponding
exponents, as in the previous case. 

 As for cases 2 and 4 of the theorem, and the component $V_2$, set $d'= p-1-(a-b)$.
Then one makes a similar computation with the monomials 
 $X^{d'/2}Y^{d'/2}$, which corresponds to $\mu_2$, and $X^{d'/2-1}Y^{d'/2+1}$, which corresponds to
 $\mu_4$.
\end{txt}

\begin{rmk} It is natural to ask whether or not the conditions in Theorem \ref{evencharcase} are
necessary as well as sufficient. This is not entirely clear, since it may well be that the 
$\chi$-component of some $v_\mu$ is divisible by a prime above $p$, which means $v_\mu$ is a test
vector in characteristic zero but not in characteristic $p$. However, it is clear
from the proof of the theorem that the conditions stated are both necessary and sufficient 
to determine test vectors in characteristic $p$. In fact, one can say a bit more that
Theorem \ref{evencharcase} in characteristic zero as well; we refer the reader to Theorem \ref{gauss-sum-intro}
in the introduction
for the statement, and the proof of that theorem below. We note also that there is some flexibility in the choice of test vector; one can
just as well use the monomial $X^{d/2+1}Y^{d/2-1}$ instead of $X^{d/2-1}Y^{d/2+1}$, which leads to a different
test vector.   
\end{rmk}

\begin{txt} We now give the proof in the odd case. Let the notation be as in Theorem \ref{oddcharcase},
and suppose that $a-b > 1$. 
We start with the component $V_1=\det^{b+1}\operatorname{Sym}^{a-b-2}(k_L^2)$. 
Set $t=(a-b-1)/2$, and consider the monomial $X^{t-1}Y^t$, corresponding to the character $\mu_1(x)
= x^{\frac{a+b-1}{2}}$. In terms of $U, V$, we get $$X^tY^{t-1}=\sum_{j=0}^{t-1} \binom{t-1}{j}
C_jU^{2j}V^{2(t-1-j)}(U+V),$$ where $C_j$ is some nonzero constant.
This expression contains $2t$ distinct monomials with 
nonzero coefficients, corresponding to distinct characters of $T$, as in the theorem. The argument for 
$V_2$ is similar, using the monomial $X^{t'-1}Y^{t'}$, where $2t'-1=p-1 - (a-b)$,
as is the case of $a-b=1$, where there is only one component.
\end{txt}

\begin{rmk} Again, there is some flexibility in the choice of test vectors: for instance one could use
the monomials $X^tY^{t-1}$ instead of $X^{t-1}Y^t$. 
\end{rmk}

\begin{rmk}
\label{onlyfp}
 A similar analysis to that carried out for the discrete series of $GL_2(\ff_p)$ can be made
 in the case of a nontrivial extension $\ff_{p^r}$. It is evident that it
 would be more complicated to carry out such an analysis in full generality,
  since the discrete series representations
of $\ff_q$ typically have many more components, and the binomial coefficients which show up may in principle
be divisible by $p$. However, it should be fairly straightforward
to deal with any specific character, for instance the trivial character of $\ff_{q^2}^*$. We have 
not pursued this question.
\end{rmk}

\subsection{Test vectors in characteristic zero: proof of the main theorems}

\begin{txt} We now return to the situation and notation
described in the introduction. Thus let $\Pi$ denote a depth
zero supercuspidal representation of $GL_2(\qq_p)$, with central character $\omega$. 
Let $\rho:GL_2(\zz_p)\rightarrow GL_{q-1}(\cc)$
denote the type. Thus $\rho$ is the inflation from $GL_2(\ff_p)$ of the $q-1$-dimensional representation
$\rho(\nu)$ where $\nu:\ff_{p^2}^*\rightarrow \cc^*$ is a character which satisfies $\nu\neq\nu^p$. 
We write $\Pi=\Pi(\nu)$. 
Let $K/\qq_p$ denote the unramified quadratic extension of $\qq_p$, and 
Let $\chi$ denote a character of $K^*$ such that $\chi$ agrees with $\omega$ on $\qq_p^*$, and assume
that $\chi$ is trivial on the principal units $U=1 + p\oo_K$. As we have seen above, $\Pi(\nu)$
contains a $p-1$-dimensional subspace $V$ stable under $GL_2(\zz_p)$ such that 
$GL_2(\zz_p)$ acts via $\rho(\nu)$, inflated from $GL_2(\ff_p)$.
We may identify $\chi$ with a character of $\ff_{p^2}^*$
as explained above, and $V$ contains a line where the action of $\ff_{q^2}^*$ viewed as a subgroup
of $GL_2(\ff_p)$ is given by the character $\chi$. 

We want to find a test vector $v\in V$, of the form $v=v_\mu$ as described in the previous paragraph.
The notation is taken to mean that $v$ is an eigenvector for matrices of the form 
$\begin{pmatrix} a & 0 \\ 0 & 1\end{pmatrix}\subset GL_2(\ff_p)$, 
with eigenvalue given by the character $a\mapsto \mu(a)$ of $\ff_p^*$. This means simply that
$v_\mu$ has nonzero projection on the $\chi$-eigenspace of the nonsplit torus $\ff_{p^2}^*$. The 
theorems of the previous section give a recipe for $\mu$ in terms of $\chi$, provided that one
has a way of labelling the characters correctly. Choose a generator $\sigma$ of the cyclic group
$\ff_{p^2}$, and identify $\sigma$ with its Teichm\"ller lift to the Witt ring $W=\oo_K$. Let
$\tau$ denote an arbitrary generator of $\mu_{p^2-1}\subset \cc^*$.  Let $\zeta=\nu(\sigma)\in\mu_{p^2-1}$.
Then we can write $\zeta=\tau^{a+ pb}$, with $0\leq a, b\leq p-1$, and at the cost
of interchanging $\nu$ and $\nu^p$, we may assume $a > b$. Further, we can write $\chi(\sigma)\in \mu_{p^2-1}$
in the form $\zeta=\tau^{r+ ps}$ with $0\leq r, s\leq p-1$. Thus we are in the situation of the previous 
paragraph, and the results given there lead immediately to the Theorems
\ref{even-case-intro} and \ref{odd-case-intro} of the Introduction. 
\end{txt}

\begin{txt} We now want to give the proof of Theorem \ref{gauss-sum-intro}. Thus assume that $\Pi$ is such
that the central character $\omega$ is trivial. We want to determine when the vectors $v_{\text{triv}}$
and $v_\text{quad}$ are test vectors for a given $\chi$ of conductor $\leq 1$, and prove the necessary condition in
the last statment of the theorem. 
Let $\xi$ denote a non-square element of 
$\ff_{p}^*$, and let $x_\xi=\begin{pmatrix} 0 & 1 \\ \xi & 0\end{pmatrix}$. Then we have 
\begin{equation}
\label{gausseval}
\begin{pmatrix} 0 & 1 \\ \xi & 0\end{pmatrix}= \begin{pmatrix} 0 & 1 \\ -1 & 0\end{pmatrix}
\begin{pmatrix} -\xi & 0 \\ 0 & 1\end{pmatrix}.
\end{equation}
We contend that the vectors $v_{\text{triv}}$ and $v_{\text{quad}}$ are eigenvectors for $x_\xi$,
with eigenvalues given by $\nu(\sqrt{\xi})=\sqrt{\xi}^{k_\nu}$ in each case. To see this,
note first of all that $v_{\text{triv}}$ and $v_{\text{quad}}$ are eigenvectors for 
$\begin{pmatrix} -\xi & 0 \\ 0 & 1\end{pmatrix}$, with eigenvalues $1$ and $-\mu_{\text{quad}}(-1)$,
respectively. We will show that each of  $v_{\text{triv}}$ and $v_{\text{quad}}$  is an eigenvector for  
$\begin{pmatrix} 0 & 1 \\ -1 & 0\end{pmatrix}$ and to calculate the corresponding eigenvalue.
This can be done by using the formula (12) on page 38 of \cite{ps-book}, which gives a formula
for the action of $w'=\begin{pmatrix} 0 & 1 \\ -1 & 0\end{pmatrix}$ under $\rho(\nu)$ as follows:
$$\left(\rho(\nu)(w')\cdot \mu\right) (y)= \sum_{x\in\ff_p^*} j(yx) \mu(x)$$
where $j$ is the Bessel function whose definition we shall recall below. Note
here that we are assuming that $\nu$ is trivial on $\ff_p^*$. Note also that
Piatetski-Shapiro is writing $\mu$ for the vector we have denoted by $v_\mu$; the character
$\mu$ itself is an eigenvector for $D^0$ in his model of $\rho(\nu)$ in the space of functions
on $\ff_p^*$. 

To spell this out, let
$\psi$ denote a nontrivial additive character of $\ff_p$. Then we have
$$j(u)=\frac{-1}{p}\sum_{\textbf{N}x=u}\psi(\tr(x))\nu(x)$$
where $u\in\ff_p^*$ and 
the sum is taken over elements $x$ of $\ff_{p^2}$ with norm $\textbf{N}x=u$. The reader
should note that the minus sign in the formula for $j$ given above is missing from the definition in
equation (4) in \cite{ps-book}; the sign makes a reappearance in the proof of the subsequent identity (6) -- see the left
side of the final formula in the proof on page 38, two lines above equation (12). The reader
could also consult \cite{bum97}, page 427, and notice the minus sign ($\epsilon = -1$) in the formula
for the action of the Weyl element given there.

Thus we get
\begin{align*}
\left(\rho(\nu)(w')\cdot \mu\right) (y)=&\frac{-1}{p}\sum_{x\in\ff_p^*} \sum_{\textbf{N}u=yx}
\psi(\tr(u))\nu(u)\mu(x)\\
= &\frac{-1}{p} \mu^{-1}(y)\sum_{x\in\ff_p^*} \sum_{\textbf{N}u=yx}
\psi(\tr(u))\nu(u)\mu(xy)\\
=& \frac{-1}{p}\mu^{-1}(y)G(\nu\cdot(\mu\circ\textbf{N}), \psi).
\end{align*}
Here $G(\nu\cdot(\mu\circ\textbf{N}), \psi)$ is the standard Gauss sum of the character 
$\nu\cdot\mu\circ\textbf{N}$
with respect to the additive character $\psi\circ\tr$. Now, if $\eta$ is any character of 
$\ff_{p^2}^*$ which satisfies $\eta=1$ on $\ff_p^*$, then we calculate that
\begin{align*}
G(\eta, \psi)=& \sum_{x\in\ff_{p^2}}\psi(\tr(u))\eta(u)\\
= & \sum_{x\in\ff_p^*\backslash\ff_{p^2}}^*\;\sum_{y\in\ff_p^*}\psi(\tr(xy))\eta(xy)\\
=& \sum_{x\in\ff_p^*\backslash\ff_{p^2}}\;\sum_{y\in\ff_p^*}\psi(y\tr(x))\eta(x)
\end{align*}
Suppose $\tr(x)\neq 0$. Then the sum on $y$ of $\psi(y\tr(x))$ has the value $-1$. 
On the other hand, if $\tr(x)=0$, then the sum on $y$ has value $p-1$. All such
$x$ with trace zero are equal to $\sqrt{\xi}$, up to multiples of $\ff_p^*$.  Thus we find
that  $$\sum_{x\in\ff_p^*\backslash\ff_{p^2}}\sum_{y\in\ff_p^*}\psi(y\tr(x))\eta(x)=
(p-1)\eta(\sqrt\xi)-\sum_{\stackrel{x\in\ff_p^*\backslash\ff_{p^2}}{\tr(x)\neq 0}}\eta(x) $$ where the final
sum is taken over the $p-1$ classes elements of nonzero trace, modulo scaling
by $\ff_p^*$. But now if we sum over all classes, then we get
$\sum_{x\in\ff_{p^2}^*/\ff_p^*}\eta(x)=0$ since $\eta$
is a nontrivial character which is trivial on $\textbf{F}_p^*$, and so we conclude
that $G(\eta, \psi) = p\eta(\sqrt\xi)$. 

Thus we obtain $$\left(\rho(\nu)(w')\cdot \mu\right) (y)= -(\nu\cdot(\mu\circ\textbf{N}))(\sqrt\xi)\mu^{-1}.$$
When $\mu^2=1$, we see that $v_\mu$ (or $\mu$, in Piatetski-Shapiro's notation) is an eigenvector
for $w'$ with eigenvalue $-(\nu\cdot\mu\circ\textbf{N})(\sqrt\xi)=\pm 1$. When $\mu$ is trivial, this eigenvalue
is $-\nu(\sqrt{\xi})$, and when $\mu$ is nontrivial quadratic, this eigenvalue is $-\nu(\sqrt{\xi})\cdot \mu\circ\textbf{N}(\sqrt{\xi})
= \nu(\sqrt{\xi})\mu(-\xi)= - \nu(\sqrt{\xi}) \mu(-1)$. 
 Going back to equation (\ref{gausseval}), we find
that the eigenvalue of $x_\xi$ on $v_{\text{triv}}$
and $v_\text{quad}$ is $\nu(\sqrt{\xi})=-\epsilon_p(\nu)$ in each case. 

In view of this computation, it is clear that the projection of $v_{\text{triv}}$
and $v_\text{quad}$ to the eigenspace of any character $\chi$ with
$\epsilon_p(\chi)=\epsilon_p(\nu)$. This completes the proof of Theorem \ref{gauss-sum-intro}.

\end{txt}
\bibliography{/Users/vvatsal/ownCloud/WorkDocs/CurrDocs/general.bib}
\bibliographystyle{amsplain}

\end{document}